\documentclass[12pt]{article}

\usepackage[english]{babel}
\usepackage[T1]{fontenc}
\usepackage{amsmath}
\usepackage{amsfonts}
\usepackage{amssymb}
\usepackage[utf8]{inputenc}
\usepackage{graphicx,color,url}
\usepackage[margin=2.5cm]{geometry}

\def\no{\noindent}
\newtheorem{lem}{Lemma}

\newtheorem{theo}{Theorem}
\newtheorem{rem}{Remark}
\def\RR{\mathbb{R}}
\def\bz{\mathbf{z}}
\def\bx{\mathbf{x}}
\def\bq{\mathbf{q}}
\def\bp{\mathbf{p}}
\def\bb{\mathbf{b}}
\def\bc{\mathbf{c}}
\def\be{\mathbf{e}}
\def\bg{\mathbf{g}}
\def\by{\mathbf{y}}
\def\bu{\mathbf{u}}
\def\bw{\mathbf{w}}

\def\tq{\tilde{q}}
\def\tp{\tilde{p}}
\def\tH{\tilde{H}}
\def\bfi{\boldsymbol{\varphi}}
\def\bffi{\boldsymbol{\phi}}
\def\bsig{\boldsymbol{\sigma}}
\def\bgam{\boldsymbol{\gamma}}
\def\bet{\boldsymbol{\eta}}
\def\calH{\mathcal{H}}
\def\calL{\mathcal{L}}
\def\calP{\mathcal{P}}
\def\calI{\mathcal{I}}
\def\Dx{\Delta x }
\def\hgam{\hat{\gamma}}
\def\pmatrix{\left(\begin{array}}
\def\endpmatrix{\end{array}\right)}
\def\dd{\mathrm{d}}
\def\ii{\mathrm{i}}

\def\dbq{\dot{\bq}}
\def\dbp{\dot{\bp}}
\def\dtq{\frac{\dd}{\dd t}{\tq}}
\def\dtp{\frac{\dd}{\dd t}{\tp}}
\def\dby{\dot\by}
\def\dbsig{\dot{\bsig}}
\def\dbu{\dot{\bu}}
\def\tJ{\tilde{J}}

\def\ddgam{\ddot{\gamma}}
\def\ddeta{\ddot{\eta}}
\def\bom{{\boldsymbol{\omega}}}
\def\proof{\underline{Proof}\quad}
\def\QED{\mbox{\,$\Box{~}$}}
\def\dgam{\dot{\gamma}}
\def\deta{\dot{\eta}}
\def\sech{\mathrm{sech}}
\def\atan{\mathrm{atan}}

\title{Energy conservation issues in the numerical solution of the semilinear wave equation}

\author{L.\,Brugnano\,$^a$\quad G.\,Frasca Caccia\,$^a$\quad F.\,Iavernaro\,$^b$\\[.5cm]
\small
$^a$\,Dipartimento di Matematica e Informatica ``U.\,Dini'', Universit\`a di Firenze, Italy\\
\small
$^b$\,Dipartimento di Matematica, Universit\`a di Bari, Italy}

\begin{document}

\maketitle

\abstract{In this paper we discuss energy conservation issues related to the  numerical solution of the semilinear wave equation.  As is well known, this problem can be cast as a Hamiltonian system that may be autonomous or not, depending on the prescribed boundary conditions.
We relate the conservation properties of the original problem to those of its semi-discrete version obtained by the method of lines. Subsequently,   we show that the very same properties can be transferred to the solutions of the fully discretized problem, obtained by using energy-conserving methods in the HBVMs (Hamiltonian Boundary Value Methods) class. Similar arguments hold true for different types of Hamiltonian Partial Differential Equations, e.g., the nonlinear Schr\"odinger equation.

\medskip
\no{\bf Keywords:} semilinear wave equation; Hamiltonian PDEs; energy-conserving methods; Hamiltonian Boundary Value Methods; HBVMs.

\medskip
\no{\bf AMS:} 65P10, 65L05, 65M20.}

\section{Introduction}\label{intro}

In this paper we discuss energy-conservation issues for the semilinear wave equation, though the approach can be extended to different kinds of {\em Hamiltonian Partial Differential Equations} (like, e.g., the nonlinear Schr\"odinger equation).  For simplicity, but without loss of generality, we shall consider the following 1D case,
\begin{eqnarray}\nonumber 
u_{tt}(x,t)&=&u_{xx}(x,t)-f'(u(x,t)), \qquad (x,t)\in (0,1)\times(0,\infty),\\ \label{wave}
u(x,0)&=&\psi_0(x),\\ \nonumber
u_t(x,0)&=&\psi_1(x), \qquad x\in(0,1),
\end{eqnarray}
coupled with suitable boundary conditions. As usual, subscripts denote  partial derivatives.
In (\ref{wave}), the functions $f$, $\psi_0$ and $\psi_1$ are supposed to be suitably regular, so that they define a regular solution $u(x,t)$ ($f'$ denotes the derivative of $f$).   The problem is completed by assigning suitable boundary conditions which we shall, at first, assume to be periodic,
\begin{equation}\label{perbc}
u(0,t)=u(1,t), \qquad t>0.
\end{equation}
In such a case, we will assume that $\psi_0$, $\psi_1$, and $f$ are such that the resulting solution also satisfies
\begin{equation}\label{uxper}
u_x(0,t)=u_x(1,t), \qquad t>0.
\end{equation}
Later on, we shall also consider the case of Dirichlet boundary conditions,
\begin{equation}\label{assbc}
u(0,t) = \varphi_0(t), \qquad u(1,t) = \varphi_1(t), \qquad t>0,
\end{equation}
and Neumann boundary conditions,
\begin{equation}\label{assnbc}
u_x(0,t) = \phi_0(t), \qquad u_x(1,t) = \phi_1(t), \qquad t>0,
\end{equation}
with $\varphi_0(t)$, $\varphi_1(t)$, $\phi_0(t)$, and $\phi_1(t)$ suitably regular.  We set 
\begin{equation}\label{vequt}
v = u_t,
\end{equation}
and define the functional
\begin{equation}\label{H1}
\calH[u,v](t)=\int_0^1 \left[\frac{1}{2}v^2(x,t)+\frac{1}{2}u_x^2(x,t)+f(u(x,t))\right]\dd x \equiv \int_0^1 E(x,t)\,\dd x.
\end{equation}
As is well known, we can rewrite (\ref{wave}) as the infinite-dimensional Hamiltonian system (for brevity, we neglect the arguments of the functions $u$ and $v$)
\begin{equation}\label{Hwave}
\bz_t=J\frac{\delta \calH}{\delta \bz}, 
\end{equation}
where
\begin{equation}\label{Hwave1}
J= \pmatrix  {cc}
0 & 1 \\
-1 & 0
\endpmatrix,
\qquad
\bz=\pmatrix{c} u\\v\endpmatrix,
\end{equation}
and 
\begin{equation}\label{Hwave2}
\frac{\delta \calH}{\delta \bz}=\left(\frac{\delta \calH}{\delta u},\frac{\delta \calH}{\delta v}\right)^\top 
\end{equation}
is the functional derivative of $\calH$. This latter is defined as follows: given  a generic functional in the form
$$\calL[q]=\int_a^b L(x,q(x),q'(x)) \dd x,$$ its functional derivative $\frac{\delta \calL}{\delta q}$ is defined by requiring that, for every function $\xi(x)$, 
$$\int_a^b\frac{\delta \calL}{\delta q}\cdot \xi\dd x \equiv  \lim_{\varepsilon\rightarrow 0} \frac{\calL[q+\varepsilon\xi]-\calL[q]}{\varepsilon} = \left. \frac{d}{d\varepsilon}\calL[q+\varepsilon\xi] \right|_{\varepsilon=0}.
$$ 
In particular, by considering a  function $\xi$ vanishing at $a$ and $b$, one obtains:
\begin{eqnarray*}
 \int_a^b\frac{\delta \calL}{\delta q}\cdot \xi\dd x &=& \left [\frac{d}{d \varepsilon}\int_a^bL(x,q+\varepsilon\xi,q'+\varepsilon \xi')\dd x \right]_{\varepsilon=0} =~
 \int_a^b\left(\frac{\partial L}{\partial q}\xi+\frac{\partial L}{\partial q'}\xi'\right)\dd x \\ &=& \int_a^b\left[ \frac{\partial L}{\partial q}\xi+ \frac{d}{dx}\left(\frac{\partial L}{\partial q'}\xi\right)-\left(\frac{d}{dx}\frac{\partial L}{\partial q'}\right)\xi\right]\dd x\\
&=& \int_a^b\left[ \frac{\partial L}{\partial q}\xi-\left(\frac{d}{dx}\frac{\partial L}{\partial q'}\right)\xi\right]\dd x~=~ \int_a^b\left( \frac{\partial L}{\partial q}-\left(\frac{d}{dx}\frac{\partial L}{\partial q'}\right)\right)\xi \,\dd x.
\end{eqnarray*}
Consequently,  
\begin{equation}\label{derfun}
\frac{\delta \calL}{\delta q}=\frac{\partial L}{\partial q}-\left(\frac{d}{dx}\frac{\partial L}{\partial q'}\right).
\end{equation}
Exploiting (\ref{derfun}), one easily verifies that  (\ref{Hwave})--(\ref{Hwave2}) are equivalent to (\ref{wave}): 
$$\bz_t=\pmatrix{c}u_t\\ v_t\endpmatrix=J\frac{\delta \calH}{\delta \bz}=\pmatrix{c}\frac{\delta \calH}{\delta v}
\\[2mm]-\frac{\delta \calH}{\delta u}\endpmatrix=\pmatrix{c}v\\ u_{xx}- f'(u)\endpmatrix,$$
or
\begin{eqnarray}\nonumber
u_t (x,t)&=& v(x,t),   \qquad (x,t)\in (0,1)\times(0,\infty),\\ \label{wave1} 
v_t(x,t)&=&u_{xx}(x,t)-f'(u(x,t)),
\end{eqnarray}
that is, the first-order formulation of the first equation in (\ref{wave}). 

The numerical treatment of Hamiltonian PDEs such as (\ref{wave}) has been the subject of an intense research activity during the past decade (see, e.g., \cite{BrRe06} for a survey). The extension of ideas and tools related to geometric integration of ordinary differential equations (ODEs) has led to the definition and analysis of various structure preserving algorithms suitable for specific or general classes of PDEs.  Two main lines of investigations are based on a multisymplectic reformulation of the equations or their semi-discretization by means of the method of lines. 

Multisymplectic structures generalize  the classical Hamiltonian structure of a Hamiltonian ODE by assigning a distinct symplectic operator for each unbounded space direction and time \cite{Br97}. A clear advantage of this approach is that it allows for an easy generalization from symplectic to multisymplectic integration.  Multisymplectic integrators are numerical methods which precisely conserve a discrete space-time symplectic structure of Hamiltonian PDEs \cite{MPS99, BrRe01, IsSc04, FMR06, Fr06, HuDeHaZh13} (a backward error analysis of such schemes may be found in \cite{MoRe03, IsSc05, IsSc06}).

In the method of lines approach, the spatial derivatives are usually approximated by finite differences or by discrete Fourier transform and the resulting system is then integrated in time by a suitable ODE integrator. Spectral methods have revealed very good potentialities especially in the case of periodic boundary conditions \cite{FW78,WMGSS91}.\footnote{They have been also applied to multisymplectic PDEs \cite{BrRe01a,ChQi01,Wa07}.}  For weakly nonlinear term $f'$ in (\ref{wave}), the modulated Fourier expansion technique \cite[Chapter XIII]{HLW06} has been adapted to both the semi-discretized and the full-discretized systems to  state long-time near conservation of energy, momentum, and actions \cite{HL08,CHL08}. In general,  quoting \cite[p.\,187]{SSCa94}, if the PDEs are of Hamiltonian type, (\dots) the space discretization should be carried out in such a way that the resulting system of ODEs is Hamiltonian (for a suitable Poisson bracket) and the time integration should also be carried out by a symplectic or Poisson integrator. This approach (which we shall consider here), has been the subject of many researches (e.g., \cite{Cano06,FrOrSS90,HeAb93,Hu91,LiQin88,LuSc97,OWW04,QinZhang90}). Whichever is the considered discretization, the main aim is that of keeping conserved discrete counterparts of continuous invariants, as done, e.g., in \cite{Fu01,FM11,KF09,Ma07,MSFM02}, with the so called {\em discrete variational derivative method}. Additional references are  \cite{Fla92,LiVQ95,StVa78}.

In this paper, we focus our attention on numerical techniques able to provide a full discretization of the original system with the discrete energy behaving consistently with the energy function associated with (\ref{wave}). More precisely, to approximate the second order spatial derivative, we use either a central finite difference or a spectral expansion, and then we derive a semi-discrete analogue of the conservation law associated with the energy density. As is well known, whatever the boundary conditions,  the rate of change of the energy density integrated over an interval depends only on the flux through its endpoints. We show that the use of an energy-conserving method to discretize the time assures a precise reproduction of the above mentioned conservation law of the semi-discrete model. In particular, if there is no net flux into or out of the interval, then the integrated energy density is precisely conserved, meaning that it remains constant over time.  Some of the presented results, in the case of periodic boundary conditions, are already known (see, e.g., \cite{QinZhang90,OWW04,Cano06}). Nevertheless, these authors mainly focus on the conservation properties of the semi-discrete model, and consider accurate symplectic integrators for their solution. Instead, we are here more interested in a precise conservation of the semi-discrete energy and, because of this reason, we consider energy-conserving methods.  Moreover, the algebraic form in which we cast the semi-discrete problem is quite concise, and allows for a simple extension to the case where the boundary conditions are not periodic. 

To the best of our knowledge, only the case of periodic boundary conditions has been studied thoroughly.  In such a case, the integral of $E$ (see (\ref{H1})) is indeed a conserved quantity  and one obtains {\em energy conservation}  (see, e.g., \cite{LeRe04}). Therefore, it makes sense to look for a corresponding conservation property, when numerically solving the problem (as done, e.g., in  \cite{QinZhang90,OWW04,Cano06}). Nevertheless, also in the other cases, which are of interest in applications, the qualitative properties of the solution can be suitably reproduced in the discrete approximation by slightly generalizing the arguments.  In fact, in all cases, one may derive a semi-discrete problem which turns out to be Hamiltonian, and whose Hamiltonian mimics a semi-discrete energy which is exactly conserved. Consequently, it makes sense to use energy-conserving methods for their numerical solution.

Energy conserving methods, in turn, have been the subject of many investigations, in the ODE setting, during the past years: we quote, as an example, {\em discrete gradient methods} \cite{McLQR99,McLQ14}, {\em time finite elements} \cite{BS00,BS00_1}, the {\em average vector field method} \cite{QMcL08,CMcLMcLOQW09,COS14} and its generalizations \cite{Ha10}. This latter method has also been considered in the PDE setting (e.g., \cite{CGMcLMcLONOQ12}). In particular, we shall here consider the energy-conserving methods in the class of {\em Hamiltonian Boundary Value Methods (HBVMs)} \cite{BIT09,BIT09_1,BIT10,BIT11,BIT12,BIT12_1,BFI13}, which are methods based on the concept of {\em discrete line integral}, as defined in \cite{IP07,IP08,IT09}. Such methods have been also generalized to the case of  different conservative problems \cite{BCMR12,BI12,BIT12_2,BIT12_3,BS13} and, more recently, they have been used for numerically solving Hamiltonian boundary value problems \cite{ABI15}.  

With this premise, the paper is organized as follows:

\begin{itemize}
\item we study, at first, the discrete problems derived by a finite-difference spatial discretization.
In particular, in Section~\ref{persec} we study the case in which problem (\ref{wave}) is completed by the periodic boundary conditions (\ref{perbc}); the case of Dirichlet boundary conditions (\ref{assbc}) will be the subject of Section~\ref{assec}; at last, the case of Neumann boundary conditions (\ref{assnbc}) will be examined in Section~\ref{assnec};

\item we then study, in Section~\ref{pbcr}, the case where a Fourier-Galerkin space discretization is considered. For sake of brevity, we shall here be concerned only with the case where periodic boundary conditions are prescribed.  Also, higher order finite-difference approximations are sketched;

\item  in Section~\ref{implement} we sketch the efficient implementation of the proposed energy-conserv\-ing methods. In Section~\ref{numtest}  we report a few numerical tests,  whereas Section~\ref{final} contains a few concluding remarks;

\item  finally, in the Appendix we sketch the way how the whole approach can be extended to different kinds of Hamiltonian PDEs. In particular, we consider the nonlinear Schr\"odinger equation.

\end{itemize}

\section{The case of periodic boundary conditions}\label{persec}
 By considering that the time derivative of the integrand function $E(x,t)$ defined at (\ref{H1}) satisfies (see (\ref{wave1}))
\begin{eqnarray*}
E_t(x,t)&=&v(x,t)v_t(x,t)+u_x(x,t)u_{xt}(x,t)+f'(u(x,t))u_t(x,t) \\ &=& v(x,t)(u_{xx}(x,t)-f'(u(x,t)))+u_x(x,t)v_x(x,t)+f'(u(x,t))v(x,t)\\
&=&v(x,t)u_{xx}(x,t)+u_x(x,t)v_x(x,t)=(u_x(x,t)v(x,t))_x ~\equiv~-F_x(x,t),\end{eqnarray*}
one derives the conservation law:
\begin{equation}\label{EtFx}
E_t(x,t)+F_x(x,t)=0,\qquad\text{with}\qquad F(x,t)=-u_x(x,t)v(x,t).
\end{equation}
Consequently, because of the periodic boundary conditions (\ref{perbc}) (and (\ref{uxper})), one obtains
%\begin{equation}\label{derH}
$$\dot{\calH}[\bz](t)=\int_0^1 E_t(x,t)\dd x=[u_x(x,t)v(x,t)]_{x=0}^1 = 0,$$
%\end{equation}
where, as usual, the dot denotes the time derivative. Therefore (\ref{H1}) is a conserved quantity, so that at $t=h$ one has:
$$\calH[\bz](h) = \calH[\bz](0).$$
We also recast the Hamiltonian function in a more convenient form to be used in the sequel. In case of the periodic boundary conditions (\ref{perbc}), from (\ref{H1}) one has 
\begin{eqnarray}\nonumber
\calH[\bz](t)&\equiv&\int_0^1 E(x,t)\,\dd x~=~\int_0^1 \left[\frac{1}{2}v^2(x,t)+\frac{1}{2}u_x^2(x,t)+f(u(x,t))\right]\dd x\\  \nonumber
&=&\int_0^1 \left[\frac{1}{2}v^2(x,t)+\frac{1}{2}[(u(x,t)u_x(x,t))_x-u(x,t)u_{xx}(x,t)]+f(u(x,t))\right]\dd x\\ \nonumber
&=&\int_0^1 \left[\frac{1}{2}v^2(x,t)-\frac{1}{2}u(x,t)u_{xx}(x,t)+f(u(x,t))\right]\dd x +  \frac{1}{2}\underbrace{\left[u(x,t)u_x(x,t)\right]_{x=0}^1}_{=0}\\
&=&\int_0^1 \left[\frac{1}{2}v^2(x,t)-\frac{1}{2}u(x,t)u_{xx}(x,t)+f(u(x,t))\right]\dd x,\label{constH}
\end{eqnarray}
where $\left[uu_x\right]_{x=0}^1=0$ because of the periodic boundary conditions (\ref{perbc}) (and (\ref{uxper})).
\subsection{Semi-discretization}\label{semiper}
For numerically solving problem (\ref{wave})-(\ref{perbc}), let us introduce the following discretization of the space variable,\begin{equation}\label{Dx1}x_i=i\Dx ,\quad  i=0,\ldots, N, \qquad \Dx =1/N,\end{equation} and the vectors: 
$$\bx =\pmatrix{c}x_0\\ \vdots\\x_{N-1}\endpmatrix, \quad\bq(t) =\pmatrix{c}u_0(t)\\ \vdots\\u_{N-1}(t)\endpmatrix,\quad\bp(t) =\pmatrix{c}v_0(t)\\ \vdots\\v_{N-1}(t)\endpmatrix\in\mathbb{R}^{N},$$ with 
\begin{equation}\label{uivi}
u_i(t)\approx u(x_i,t), \qquad v_i(t)\approx v(x_i,t)\equiv u_t(x_i,t).
\end{equation}
Because of the periodic boundary conditions (\ref{perbc}), we also set:
%\begin{equation}\label{cc}
$$u_{N}(t)\equiv u_0(t),\quad u_{-1}(t)\equiv u_{N-1}(t), \qquad t\ge0.$$
%\end{equation}
Approximating the second derivative in (\ref{wave1}) as
\begin{equation}\label{uxx}
u_{xx}(x_i,t)\approx \frac{u_{i+1}(t)-2u_i(t)+u_{i-1}(t)}{\Dx ^2},\qquad i=0,\ldots, N-1,
\end{equation}
yields the following semi-discrete problem
\begin{eqnarray}\label{wave1d}
\dbq &=& \bp,\\
\dbp &=& -\frac{1}{\Dx ^2} T_N \bq- f'(\bq ), \qquad t>0,\nonumber
\end{eqnarray}
with the initial condition
\begin{equation}\label{q0p0}
\bq(0) = \psi_0(\bx), \qquad \bp(0)=\psi_1(\bx),
\end{equation}
(with an obvious meaning for $f'(\bq)$, $\psi_0(\bx)$, and $\psi_1(\bx)$)
and the following approximation of the Hamiltonian (\ref{constH}),
\begin{equation}\label{Hperd}
H\equiv H(\bq,\bp) =\Dx \left[\frac{\bp ^\top \bp}2 +\frac{\bq ^\top T_N\bq}{2\Dx ^2} +\be ^\top f(\bq )\right],
\end{equation}
where $T_N$ is a circulant matrix,\footnote{Because of the periodic boundary conditions (\ref{perbc}).}
\begin{equation}\label{TNc}
T_N=\left[\begin{array}{ccccc}
2 & -1 &   &  & -1\\
-1 & \ddots & \ddots  & & \\
 & \ddots & \ddots & \ddots &  \\
 & & \ddots& \ddots& -1\\
-1 & & & -1&  2\\
\end{array}\right]\in\mathbb{R}^{N\times N},
\end{equation}
and
\begin{equation}\label{be}
\be=\pmatrix{ccc}1& \dots &1\endpmatrix^\top \in\mathbb{R}^N.
\end{equation}
Problem (\ref{wave1d}) is clearly Hamiltonian. In fact, one has
%\begin{equation}\label{wave1dH}
$$\dot\bq = \frac{1}{\Dx } \nabla_\bp H,\qquad \dot\bp = -\frac{1}{\Dx } \nabla_\bq H,$$
%\end{equation}
or, by introducing the vector
%\begin{equation}\label{by}
$$\by =\pmatrix{c}\bq\\ \bp\endpmatrix,$$
%\end{equation}
one obtains the more compact form
\begin{equation}\label{wave1dy}
\dby  = J_N \nabla H(\by), \qquad\mbox{with}\qquad J_N=\frac{1}{\Dx}\pmatrix{cc} &I_N\\ -I_N\endpmatrix,
\end{equation}
where here and in the sequel we use,  when appropriate, the notation $H(\by)=H(\bq,\bp)$. Consequently,
$$\dot{H}(\by) = \nabla H(\by)^\top \dby =  \nabla H(\by)^\top J_N\nabla H(\by) =0,$$ 
because $J_N$ is skew-symmetric. One then concludes that the discrete approximation (\ref{Hperd}) to (\ref{constH}) is a conserved quantity for the semi-discrete problem (\ref{wave1dy}). Writing (\ref{Hperd})  in componentwise form, 
\begin{equation}\label{Hdisp}
H(\bq,\bp) = \Dx\sum_{i=0}^{N-1} \left(\frac{1}{2}v_i^2 - u_i\frac{u_{i-1}-2u_i+u_{i+1}}{2\Dx^2} +f(u_i) \right),
\end{equation}
one notices that  (\ref{Hperd}) is nothing but the approximation of (\ref{constH})  via the composite trapezoidal rule (provided that the second derivative $u_{xx}$ has been previously approximated as indicated at (\ref{uxx}), and taking into account the periodic boundary conditions (\ref{perbc})).  Consequently, one sees that (\ref{Hdisp}) is a $O(\Dx^2)$ approximation to (\ref{constH}).

\subsection{Full discretization}\label{disper}
Problem (\ref{wave1dy}) can be discretized by using a HBVM$(k,s)$ method which allows for an (at least {\em practical}\/) conservation of (\ref{Hperd}), by using a suitably large value $k\ge s$ \cite{BIT12_1},  as is shown in the sequel. Let us study the approximation to the solution over the time interval $[0,h]$, representing the very first step of the numerical approximation, to be repeated subsequently. For this purpose, we shall consider the orthonormal polynomial basis over the interval [0,1], $\{P_j\}$, given by the shifted and scaled Legendre polynomials: $$\deg P_i=i,\qquad\int_0^1P_i(x)P_j(x)\dd x=\delta_{ij},\qquad \forall i,j\geq 0,$$
$\delta_{ij}$ being the Kronecker symbol. Let us then expand the right-hand side of (\ref{wave1dy}) along this basis, thus obtaining
\begin{equation}\label{ydot}
\dby(ch)=\sum_{j\geq 0}\gamma_j(\by )P_j(c), \quad c\in[0,1],
\end{equation}
with
\begin{equation}\label{gammaj}
\gamma_j(\by)=\int_0^1 J_N\nabla H(\by (\tau h))P_j(\tau)\dd\tau, \quad j\geq 0.
\end{equation}
It is possible to prove the following result \cite{BIT12_1}.
\begin{lem}\label{lemma1} Assume $\nabla H(\by(\cdot))$ can be expanded in Taylor series at $0$. Then: $$\gamma_j(\by) = O(h^j) \in \RR^{2N}, \qquad j=0,1,\dots.$$\end{lem}
Setting the initial condition (see (\ref{wave}))
\begin{equation}\label{y0}
\by _0=\pmatrix{c}\psi_0(\bx )\\ \psi_1(\bx )\endpmatrix,
\end{equation} 
with $\psi_j(\bx)$, $j=0,1$, the vector whose entries are given by $\psi_j(x_i)$, the solution of (\ref{ydot})-(\ref{y0}) is then formally  given by:
\begin{equation}\label{ych}
\by (ch)=\by _0+h\sum_{j\geq 0}\gamma_j(\by )\int_0^cP_j(x)\dd x,\qquad c\in[0,1].
\end{equation}
In order to obtain a polynomial approximation $\bsig\in\Pi_s$ to (\ref{ych}), we consider the following {\em truncated} initial value problem \cite{BIT12_1},
\begin{equation}\label{dotsig}
\dbsig(ch)=\sum_{j=0}^{s-1}\gamma_j(\mathbf{\bsig })P_j(c),\qquad c\in[0,1],\qquad \bsig (0)=\by _0,
\end{equation}
where $\gamma_j(\bsig )$ is still given by (\ref{gammaj}) by replacing $\by$ with $\bsig$. The polynomial approximation to (\ref{ych}) is then formally given by:
$$\bsig (ch)=\by _0+h\sum_{j=0}^{s-1}\gamma_j(\bsig )\int_0^c P_j(x)\dd x,\qquad c\in[0,1].$$
The use of a quadrature formula of order $2k$ to approximate the integral defining $\gamma_j(\bsig)$ (see  (\ref{gammaj})) would give \cite{BIT12_1}
\begin{eqnarray}\nonumber
\gamma_j(\bsig) &=&\int_0^1J_N\nabla H(\bsig(\tau h))P_j(\tau)\dd\tau \\
&=&\underbrace{\sum_{\ell=1}^kb_\ell P_j(c_\ell)J_N\nabla H(\bsig(c_\ell h))}_{=\hgam_j(\bsig)} + \Delta_j(h)  \label{hgamj}
~\equiv~ \hgam_j(\bsig)+ \Delta_j(h),
\end{eqnarray}
with
\begin{equation}\label{Deltaj}
\Delta_j(h) = O(h^{2k-j}) \in\RR^{2N}, \qquad j=0,\dots,s-1.
\end{equation}
In such a case, however, we have a different polynomial $\bu\in\Pi_s$, in place of $\bsig$, solution of the problem 
\begin{eqnarray}\label{dotu}  
\dbu(ch)&=&\sum_{j=0}^{s-1} \hgam_j(\bu) P_j(c) ,\qquad c\in[0,1],\qquad
\bu (0)=\by _0,\\ \nonumber
\hgam_j(\bu) & =& \sum_{i=1}^k b_i P_j(c_i) J_N \nabla H(\bu(c_ih)), \qquad j=0,\dots,s-1,
\end{eqnarray}
instead of (\ref{dotsig}):  this latter problem defines a HBVM$(k,s)$ method.

If $H(\bq,\bp)$ in (\ref{Hperd}) is a polynomial of degree $\nu\ge2$ (which means that $f\in\Pi_\nu$)\footnote{Indeed,  $H$ contains at least a quadratic term.}, and $k$ is an integer such that\,
\begin{equation}\label{nusk}
k\ge \frac{1}2 \nu s \qquad \Leftrightarrow \qquad \nu\leq \frac{2k}s,
\end{equation}
we can exactly compute the integrals $\gamma_j(\bsig)$ by means of a Gauss-quadrature formula of order $2k$, so that $\bu\equiv\bsig$ and, then:
\begin{eqnarray}\label{Hsig}
\lefteqn{H(\bsig (h))-H(\bsig (0))~=~h\int_0^1\nabla H(\bsig (\tau h))^\top \dbsig(\tau h)\dd\tau}\\
\nonumber
&=&h\int_0^1\nabla H(\bsig (\tau h))^\top \sum_{j=0}^{s-1}P_j(\tau)\gamma_j(\bsig )\dd\tau
~=~h\Dx^2  \sum_{j=0}^{s-1}\gamma_j(\bsig )^\top J_N\gamma_j(\bsig )~=~0,
\end{eqnarray}
due to the fact that $J_N$ is skew-symmetric.
If $f$, and then $H$, is not a polynomial,  by taking into account (\ref{wave1dy}) and (\ref{hgamj})--(\ref{dotu}), the error on the Hamiltonian $H$, at $t=h$, is:
\begin{eqnarray}\nonumber
\lefteqn{H(\bu(h))-H(\bu(0))~=~h\int_0^1\nabla H(\bu (\tau h))^\top \dbu(\tau h)\dd\tau}\\
\nonumber
&=&h\int_0^1\nabla H(\bu (\tau h))^\top \sum_{j=0}^{s-1}P_j(\tau)\left(\gamma_j(\bu )-\Delta_j(h)\right)\dd\tau\\
\nonumber
&=&h\Dx^2\sum_{j=0}^{s-1}\left[\overbrace{\gamma_j(\bu)^\top J_N\gamma_j(\bu)}^{=0}-\gamma_j(\bu )^\top J_N\Delta_j(h)\right]\\
&=&h\underbrace{\Dx\cdot N}_{=1} \cdot O(h^{2k}) ~\equiv~ O\left(h^{2k+1}\right), \label{errHk}
\end{eqnarray}
where the last equality follows from Lemma~\ref{lemma1} and (\ref{Deltaj}). Consequently, choosing $k$ large enough allows us to approximate the Hamiltonian $H$ within full machine accuracy. Summing up all the previous arguments and taking into account the results in \cite{BIT12_1}, the following results can be proved.

\begin{theo}\label{hbvmRK} 
The HBVM$(k,s)$ method (\ref{dotu}) is the $k$-stage Runge-Kutta method with tableau
\begin{eqnarray}\label{RKform}
\begin{array}{c|c}
\bc & \calI \calP^\top\Omega\\
\hline
       & \bb^\top
       \end{array} ~& with & \left\{\begin{array}{rcl}\bb &=& \pmatrix{ccc} b_1&\dots&b_k\endpmatrix^\top\\
       \bc &=& \pmatrix{ccc} c_1&\dots&c_k\endpmatrix^\top\end{array}\right.,
       \quad \Omega = \pmatrix{ccc} b_1\\ &\ddots \\ &&b_k\endpmatrix,\\[1mm]                           
 \nonumber                                      
        &and&  \calP = \pmatrix{c} P_{j-1}(c_i)\endpmatrix,~
                                       \calI = \pmatrix{c} \int_0^{c_i}P_{j-1}(x)\dd x\endpmatrix \in\RR^{k\times s}.
\end{eqnarray}
\end{theo}

\begin{theo}\label{hbvmpbc} Assume $k\ge s$, and define ~$\by_1=\bu(h)$~ as the new approximation to ~$\by(h)$~  provided by a HBVM$(k,s)$ method used with stepsize $h$. One then obtains:
$$\by_1-\by(h) = O(h^{2s+1}),$$ that is the method has order $2s$. Moreover, with reference to (\ref{nusk}), and assuming that $f$ is suitably regular:
$$H(\by_1)-H(\by_0) = \left\{\begin{array}{ccl}
0, &\quad&\mbox{if ~$f\in\Pi_\nu$~ and ~$\nu\leq 2k/s$,}\\[.5cm]
O(h^{2k+1}), &&\mbox{otherwise.}
\end{array}\right.$$
\end{theo}

\begin{rem}\label{rem1} From this result, it follows that one can always obtain the conservation of the discrete Hamiltonian (\ref{Hperd}) when $f$ is a polynomial, by choosing $k$ large enough. Moreover, as (\ref{errHk}) suggests, also in the non-polynomial case, a {\em practical} conservation of (\ref{Hperd}) can be gained by choosing $k$ large enough, so that the approximation is within round-off errors.
As we shall see in Section~\ref{implement}, this is not a severe drawback, since the discrete problem generated by a HBVM$(k,s)$ method has dimension $s$, independently of $k$ (see also \cite{BIT11,BIT12_1,BFI13}).
\end{rem}

\section{The case of Dirichlet boundary conditions}\label{assec}

Let us now consider the case when the considered problem is given by (\ref{wave}) with the boundary conditions (\ref{assbc}).
By repeating similar steps as done in (\ref{constH}), one obtains:
\begin{eqnarray}\nonumber
\calH[\bz](t)&=&\int_0^1 E(x,t)\dd x~\equiv~\int_0^1 \left[\frac{1}{2}v(x,t)^2+\frac{1}{2}u_x(x,t)^2+f(u(x,t))\right]\dd x\\  \nonumber
&=&\int_0^1 \left[\frac{1}{2}v(x,t)^2+\frac{1}{2}[(u(x,t)u_x(x,t))_x-u(x,t)u_{xx}(x,t)]+f(u(x,t))\right]\dd x\\ \nonumber
&=&\int_0^1 \left[\frac{1}{2}v(x,t)^2-\frac{1}{2}u(x,t)u_{xx}(x,t)+f(u(x,t))\right]\dd x +  \frac{1}{2}\left[u(x,t)u_x(x,t)\right]_{x=0}^1\\ \nonumber
&=&\int_0^1 \left[\frac{1}{2}v(x,t)^2-\frac{1}{2}u(x,t)u_{xx}(x,t)+f(u(x,t))\right]\dd x  +\\ &&
\qquad \frac{1}{2}\left[ u(1,t)u_x(1,t)-u(0,t)u_x(0,t)\right].\label{newH}
\end{eqnarray}
Moreover, $\calH[\bz]$ is no more conserved because formally (\ref{EtFx}) still holds true and, then, one obtains (see also (\ref{H1}), and taking into account the boundary conditions (\ref{assbc})):
\begin{equation}\label{varia}
\dot{\calH}[\bz](t)=\int_0^1 E_t(x,t)\dd x=[u_x(x,t)v(x,t)]_{x=0}^1 = u_x(1,t)\varphi_1'(t)-u_x(0,t)\varphi_0'(t).
\end{equation}
Equation (\ref{varia}) may be interpreted as the instant variation of the energy  which is released or gained by the system at time $t$. Thus, the continuous Hamiltonian (\ref{H1}), though no more conserved, has a {\em prescribed variation in time}. From (\ref{varia}), at $t=h$ one easily obtains:
\begin{equation}\label{incH}
\calH[\bz](h)-\calH[\bz](0) = \int_0^h \dot{\calH}[\bz](t)\dd t = \int_0^h \left[u_x(1,t)\varphi_1'(t)-u_x(0,t)\varphi_0'(t)\right]\dd t.
\end{equation}

\subsection{Semi-discretization}\label{semiass}
In order for numerically solving problem (\ref{wave})--(\ref{assbc}), let us introduce the following discretization of the space variable,
\begin{equation}\label{timed}
x_i=i\Dx ,\quad  i=0,\ldots, N+1 \qquad \Dx =1/(N+1),\end{equation} and the vectors: 
\begin{equation}\label{xqp}
\bx =\pmatrix{c}x_1\\ \vdots\\x_N\endpmatrix, \quad\bq(t) =\pmatrix{c}u_1(t)\\ \vdots\\u_N(t)\endpmatrix,\quad\bp(t) =\pmatrix{c}v_1(t)\\ \vdots\\v_N(t)\endpmatrix\in\mathbb{R}^{N},
\end{equation}
with $u_i(t)$ and $v_i(t)$ formally defined as in (\ref{uivi}).
Approximating the second derivatives in (\ref{wave1}) as follows,
\begin{equation}\label{uxx1}
u_{xx}(x_i,t)\approx \frac{u_{i+1}(t)-2u_i(t)+u_{i-1}(t)}{\Dx ^2},\qquad i=1,\ldots, N,
\end{equation}
and, moreover,
\begin{equation}\label{ux1}
 u_x(1,t)\approx \frac{u_{N+1}(t)-u_N(t)}{\Dx },\qquad u_x(0,t)\approx \frac{u_1(t)-u_0(t)}{\Dx },
\end{equation}
we then arrive at the following semi-discrete version of (\ref{newH}):
\begin{eqnarray}\nonumber
H &=&\Dx \sum_{i=1}^N\left(\frac{1}{2}v_i^2-u_i\frac{u_{i-1}-2u_i+u_{i+1}}{2\Dx ^2}+f(u_i)\right)\\
&&+\frac{1}2\left[u_{N+1}\frac{u_{N+1}-u_N}{\Dx }-u_0\frac{u_1-u_0}{\Dx }\right].\label{semidisH}
\end{eqnarray}
Moreover, because of the boundary conditions (\ref{assbc}), one has:
\begin{equation}\label{cc2}
u_{0}(t)=\varphi_0(t),\qquad u_{N+1}(t)=\varphi_1(t),
\end{equation}
so that  we obtain the following semi-discrete approximation to the Hamiltonian (\ref{newH}):
$$ 
H =\Dx \sum_{i=1}^N\left(\frac{1}{2}v_i^2-u_i\frac{u_{i-1}-2u_i+u_{i+1}}{2\Dx ^2}+f(u_i)\right)+\varphi_1\frac{\varphi_1-u_N}{2\Dx }+\varphi_0\frac{\varphi_0-u_1}{2\Dx }.
$$ 
$H$ can be rewritten in vector form as 
\begin{equation}\label{Hassdv}
H \equiv H(\bq,\bp,t) = {\Dx }\left[\frac{\bp ^\top \bp }{2}+\frac{\bq ^\top T_N\bq }{2 \Dx ^2}+\be^\top f(\bq )\right]+\frac{\bfi(t)^\top \bfi(t) }{2\Dx }-\frac{\bq^\top \bfi(t) }{\Dx },
\end{equation}
where $\be $ has been defined in (\ref{be}) and, moreover:
\begin{equation}\label{TND}
T_N=\pmatrix{ccccc}
2 & -1 &   &  & \\
-1 & \ddots & \ddots  & & \\
 & \ddots & \ddots & \ddots &  \\
 & & \ddots& \ddots& -1\\
 & & & -1&  2\\
\endpmatrix\in\mathbb{R}^{N\times N},
\qquad
\bfi(t) =\pmatrix{c}
\varphi_0(t)\\
0\\
\vdots\\
0\\
\varphi_1(t)
\endpmatrix\in\mathbb{R}^N.
\end{equation}
With reference to (\ref{Hassdv})-(\ref{TND}), the corresponding semi-discrete problem is then given by:
\begin{eqnarray}\label{wave1da}
\dbq &=& \bp ~\equiv~\frac{1}\Dx\nabla_\bp H,  \qquad t>0,\\
\dbp &=& -\frac{1}{\Dx ^2} T_N \bq+\frac{1}{\Dx^2}\bfi- f'(\bq ) ~\equiv~-\frac{1}\Dx\nabla_\bq H,\nonumber
\end{eqnarray}
which is clearly Hamiltonian, though the Hamiltonian (\ref{Hassdv}) is now non-autonomous, because of the boundary conditions (\ref{assbc}). In order to conveniently handle this problem, we at first
transform (\ref{wave1da}) into an enlarged {\em autonomous} Hamiltonian system, by introducing the following auxiliary conjugate  scalar variables,
\begin{equation}\label{auxvar}
\tq \equiv t, \qquad \tp,
\end{equation}
and the augmented Hamiltonian (compare with (\ref{Hassdv})),
\begin{eqnarray}\nonumber
\tH(\bq,\bp,\tq,\tp) &=&{\Dx }\left[\frac{\bp ^\top \bp }{2}+\frac{\bq ^\top T_N\bq }{2 \Dx ^2}+\be ^\top f(\bq )\right]+\frac{\bfi (\tq)^\top \bfi(\tq) }{2\Dx }-\frac{\bq ^\top \bfi(\tq) }{\Dx } +\tp\\  \label{Hmod}
&\equiv& H(\bq,\bp,\tq)+\tp.
\end{eqnarray}
The dynamical system corresponding to this new Hamiltonian function is, for $t>0$:
\begin{eqnarray}\nonumber
\dbq &=& \bp ~\equiv~\frac{1}\Dx\nabla_\bp \tH,   \\ \nonumber
\dbp &=& -\frac{1}{\Dx ^2} T_N \bq+\frac{1}{\Dx^2}\bfi- f'(\bq ) ~\equiv~-\frac{1}\Dx\nabla_\bq \tH,\\ \label{wave1m}
\dtq  &=& 1 ~\equiv~\frac{\partial}{\partial\tp} \tH,\\ \nonumber
\dtp  &=& -\frac{\varphi_0(\tq)-u_1}{\Dx}\varphi_0'(\tq) -\frac{\varphi_1(\tq)-u_N}{\Dx}\varphi_1'(\tq)~\equiv~-\frac{\partial}{\partial\tq} \tH,
\end{eqnarray}
with initial conditions given by (see (\ref{xqp}))
\begin{equation}\label{init1}
\bq(0)=\psi_0(\bx), \qquad \bp(0)=\psi_1(\bx), \qquad \tq(0)=\tp(0)=0.
\end{equation}
The first 3 equations in (\ref{wave1m}) exactly coincides with (\ref{wave1da}) (considering that $\tq\equiv t$), whereas the last one allows for the conservation of $\tH$:
$$\tH(\bq(t),\bp(t),\tq(t),\tp(t)) = \tH(\bq(0),\bp(0),0,0) \equiv H(\bq(0),\bp(0),0),\qquad t\ge0.$$
Indeed, one readily sees that
\begin{equation}\label{dHt}
\frac{\dd}{\dd t}\tH(\bq,\bp,\tq,\tp)=\overbrace{\nabla_\bq \tH^\top \dbq + \nabla_\bp \tH^\top \dbp}^{=0} + \underbrace{\frac{\partial}{\partial \tq}\tH \dtq+ \frac{\partial}{\partial \tp}\tH \dtp}_{=0} =0,
\end{equation}
by virtue of (\ref{wave1m}). Consequently, by recalling that $\tq\equiv t$, from (\ref{Hassdv}) and (\ref{dHt}) one obtains:
$$\frac{\dd}{\dd t} H(\bq,\bp,t) = \frac{\partial}{\partial t} H(\bq,\bp,t) =\left[\frac{\varphi_0(t)-u_1}{\Dx}\varphi_0'(t) +  \frac{\varphi_1(t)-u_N}{\Dx}\varphi_1'(t)\right],$$ 
which is the discrete counterpart of (\ref{varia}), via the approximation (\ref{ux1}) and taking into account the boundary conditions (\ref{cc2}). Consequently,
one obtains the following  semi-discrete analogue of (\ref{incH}):
\begin{eqnarray}\nonumber
\lefteqn{H(\bq(h),\bp(h),h)-H(\bq(0),\bp(0),0) =}\\
&=&\int_0^h \left[\frac{u_{N+1}(t)-u_N(t)}{\Dx}\varphi_1'(t) -
\frac{u_1(t)-u_0(t)}{\Dx}\varphi_0'(t)\right]\dd t.\label{disan}
\end{eqnarray}
 
\begin{rem}\label{prescri}
 It is  clear that (\ref{disan}) is equivalent to keep constant $\tH(\bq(t),\bp(t),t,\tp(t))$ along the solution of (\ref{wave1m}).
 Consequently, by conserving the augmented Hamiltonian $\tH$, one obtains that $H$ satisfies a prescribed variation in time which, in turn, is consistent with the corresponding continuous one.
 \end{rem}

In order to simplify the notation, let us set
\begin{equation}\label{by1}
\by = \pmatrix{c} \bq\\ \bp\\ \tq\\ \tp,\endpmatrix, \qquad \tJ_N = \pmatrix{cc|cc} & \frac{1}{\Dx} I_N\\
-\frac{1}{\Dx} I_N &\\ \hline &&&1\\&&-1&\endpmatrix,
\end{equation}
so that (\ref{wave1m})-(\ref{init1}) can be rewritten as
\begin{equation}\label{wave2m}
\dby = \tJ_N\nabla\tH(\by), \quad t>0, \qquad \by(0) = (\psi_0(\bx)^\top ,\psi_1(\bx)^\top ,0,0)^\top .
\end{equation}

\subsection{Full discretization}\label{disass}

The full discretization of (\ref{by1})-(\ref{wave2m}) follows similar steps as those seen in Section~\ref{disper} for (\ref{wave1dy}). Let us then expand the right-hand side in (\ref{wave2m}) as done in (\ref{ydot})-(\ref{gammaj}), and consider the polynomial approximation of degree $s$ given by (\ref{dotsig}), by formally replacing $H$ with $\tH$. In such a case, one obtains energy conservation, since (compare with (\ref{Hsig}))
\begin{eqnarray}\nonumber
\lefteqn{\tH(\bsig (h))-\tH(\bsig (0))~=~h\int_0^1\nabla \tH(\bsig (\tau h))^\top \dbsig(\tau h)\dd\tau}\\
&=&h\int_0^1\nabla \tH(\bsig (\tau h))^\top \sum_{j=0}^{s-1}P_j(\tau)\gamma_j(\bsig )\dd\tau
\label{Hsig1}
=h  \sum_{j=0}^{s-1}\gamma_j(\bsig )^\top \tJ_N^{-\top}\gamma_j(\bsig )~=~0,
\end{eqnarray}
since
\begin{equation}\label{tJNtop}
\tJ_N^{-\top} = \pmatrix{cc|cc} & \Dx I_N\\
-\Dx I_N &\\ \hline &&&1\\&&-1&\endpmatrix
\end{equation}
is skew-symmetric. Consequently, if one is able to exactly compute the integrals, by means of a quadrature rule based at $k\ge s$ Gaussian points, with $k$ large enough, energy conservation is gained. This is the case, provided that $\tH$ is a polynomial, that is, $f\in\Pi_\nu$ and $\varphi_0,\varphi_1\in\Pi_\rho$, and, moreover, $k$ satisfies:
\begin{equation}\label{kgenuro}
k\ge \frac{1}2\max\left\{ \nu s, 2\rho+s-1, \rho+2s-1\right\}   
\end{equation}
(we observe that, in case $\rho=0$, such a bound reduces to the bound (\ref{nusk}), obtained in the case of periodic boundary conditions).
Differently, by approximating the integrals by means of a Gaussian quadrature of order $2k$, one obtains, with arguments similar to those used in  (\ref{hgamj})-(\ref{Deltaj}), 
\begin{eqnarray}\nonumber
\lefteqn{\gamma_j(\bsig) ~=~\int_0^1\tJ_N\nabla \tH(\bsig(\tau h))P_j(\tau)\dd\tau }\\
&=&\underbrace{\sum_{\ell=1}^kb_\ell P_j(c_\ell)\tJ_N\nabla \tH(\bsig(c_\ell h))}_{=\hgam_j(\bsig)} + \Delta_j(h)  \label{hgamj1}
~\equiv~ \hgam_j(\bsig)+ \Delta_j(h),
\end{eqnarray}
with
\begin{equation}\label{Deltaj1}
\Delta_j(h) = O(h^{2k-j}) \in\RR^{2N+2}, \qquad j=0,\dots,s-1.
\end{equation}
In such a case, we have again a different polynomial $\bu\in\Pi_s$, in place of $\bsig$, solution of a problem formally  still given by (\ref{dotu}) with $H$ replaced by $\tH$.
As a consequence, by taking into account (\ref{Deltaj1}), the error in the Hamiltonian $\tH$, at $t=h$, is given by (see (\ref{by1})):
\begin{eqnarray}\nonumber
\lefteqn{\tH(\bu(h))-\tH(\bu(0))}\\ \nonumber
&=&h\int_0^1\nabla \tH(\bu (\tau h))^\top \dbu(\tau h)\dd\tau
=h\int_0^1\nabla \tH(\bu (\tau h))^\top \sum_{j=0}^{s-1}P_j(\tau)\left(\gamma_j(\bu )-\Delta_j(h)\right)\dd\tau\\
\nonumber
&=&h\sum_{j=0}^{s-1}\left[\overbrace{\gamma_j(\bu)^\top \tJ_N^{-\top}\gamma_j(\bu)}^{=0}-\gamma_j(\bu )^\top \tJ_N^{-\top}\Delta_j(h)\right]=h\underbrace{\Dx\cdot N}_{<\, 1} \cdot O(h^{2k}) \equiv O\left(h^{2k+1}\right), \label{errHk1}
\end{eqnarray}
where the last equality follows from (\ref{timed}), (\ref{tJNtop}), and Lemma~\ref{lemma1}. Consequently, choosing $k$ large enough allows us to approximate the Hamiltonian $\tH$ within full machine accuracy. 

All the above arguments can be summarized by the following theorem, which generalizes Theorem~\ref{hbvmpbc}  to the present case.

\begin{theo}\label{hbvmpbc1} Assume $k\ge s$, and define ~$\by_1=\bu(h)$~ as the new approximation to $\by(h)$, solution of (\ref{by1})-(\ref{wave2m}),  provided by a HBVM$(k,s)$ method used with stepsize $h$. One then obtains:
$$\by_1-\by(h) = O(h^{2s+1}),$$ that is the method has order $2s$. Moreover, assuming that $f$, $\varphi_0$, and $\varphi_1$ in (\ref{wave})-(\ref{assbc}) are suitably regular:
$$\tH(\by_1)-\tH(\by_0) = \left\{\begin{array}{ccl}
0, &\quad&\mbox{if ~$f\in\Pi_\nu$, $\varphi_0,\varphi_1\in\Pi_\rho,$~ and ~(\ref{kgenuro}) holds true,}\\[.5cm]
O(h^{2k+1}), &&\mbox{otherwise.}
\end{array}\right.$$
\end{theo}
Clearly,  considerations similar to those stated in Remark~\ref{rem1} can be repeated also in the present situation.

\section{The case of Neumann boundary conditions}\label{assnec}
As done in the case of Dirichlet boundary conditions, also when Neumann boundary conditions are prescribed, one starts from the formulation (\ref{newH}) of the continuous Hamiltonian function, and then considers its semi-discretization  (\ref{semidisH}). In so doing, one arrives at the very same formulation (\ref{Hassdv}), with $T_N$ defined as in (\ref{TND}), whereas, by considering the Neumann boundary conditions (\ref{assnbc}), and the approximations (\ref{ux1}) used to derive (\ref{semidisH}), $\bfi(t)$ is now formally defined as follows:
$$ \bfi(t) = \pmatrix{ccccc} u_1-\phi_0(t)\Dx, &0,&\dots~,&0,& u_N+\phi_1(t)\Dx\endpmatrix^\top.$$ 
In fact, this is equivalent to use the following definitions for $u_0(t)$ and $u_{N+1}(t)$, 
\begin{equation}\label{u0Npiu1}
u_0(t) = u_1-\phi_0(t)\Dx, \qquad u_{N+1}(t) = u_N+\phi_1(t)\Dx,
\end{equation}
which we shall use in the sequel.
We prefer, however, to derive the semi-discrete Hamiltonian by following a slightly different route, as described below. In more details,
starting from (\ref{semidisH}),  one obtains:
\begin{eqnarray}\nonumber 
\lefteqn{H ~=~}\\ \nonumber
&=&  \Dx \sum_{i=1}^N\left(\frac{1}{2}v_i^2-u_i\frac{u_{i-1}-2u_i+u_{i+1}}{2\Dx ^2}+f(u_i)\right)
+\frac{1}2\left[u_{N+1}\frac{u_{N+1}-u_N}{\Dx }-u_0\frac{u_1-u_0}{\Dx }\right] \\
\nonumber
&=&\Dx \sum_{i=2}^{N-1}\left(\frac{1}{2}v_i^2-u_i\frac{u_{i-1}-2u_i+u_{i+1}}{2\Dx ^2}+f(u_i)\right)+\frac{1}2\left[u_{N+1}\frac{u_{N+1}-u_N}{\Dx }-u_0\frac{u_1-u_0}{\Dx }\right]\\
\nonumber
&&+\Dx \left(\frac{1}{2}v_1^2-u_1\frac{u_0-2u_1+u_2}{2\Dx ^2}+f(u_1) ~+~\frac{1}{2}v_N^2-u_N\frac{u_{N-1}-2u_N+u_{N+1}}{2\Dx ^2}+f(u_N)\right)\\
\nonumber
&=&\Dx \sum_{i=2}^{N-1}\left(\frac{1}{2}v_i^2-u_i\frac{u_{i-1}-2u_i+u_{i+1}}{2\Dx ^2}+f(u_i)\right)+\frac{1}2\left[\frac{(u_{N+1}-u_N)^2}{\Dx }+\frac{(u_1-u_0)^2}{\Dx }\right]\\
\nonumber
&&+\Dx \left(\frac{1}{2}v_1^2-u_1\frac{-u_1+u_2}{2\Dx ^2}+f(u_1) ~+~\frac{1}{2}v_N^2-u_N\frac{u_{N-1}-u_N}{2\Dx ^2}+f(u_N)\right),\\
\label{semidisH1}
\end{eqnarray}
which can be cast in vector form as
\begin{equation}\label{Hassnv}
H \equiv H(\bq,\bp,t) = {\Dx }\left[\frac{\bp ^\top \bp }{2}+\frac{\bq ^\top T_N\bq }{2 \Dx ^2}+\be^\top f(\bq )\right]+\frac{\bw(\bq,t)^\top \bw(\bq,t) }{2\Dx },
\end{equation}
where $\be $ has been defined in (\ref{be}), $\bq$ and $\bp$ are defined at (\ref{xqp}), whereas:
\begin{equation}\label{TNN}
T_N=\pmatrix{ccccc}
1 & -1 &   &  & \\
-1 & 2 & \ddots  & & \\
 & \ddots & \ddots & \ddots &  \\
 & & \ddots& 2& -1\\
 & & & -1&  1\\
\endpmatrix\in\mathbb{R}^{N\times N},\qquad
\bw(\bq,t) = \pmatrix{c} 
u_1-u_0(t)\\
0\\
\vdots
\\
0\\
u_{N+1}(t)-u_N
\endpmatrix\in\RR^N.
\end{equation}
We emphasize that $u_0(t)$ and $u_{N+1}(t)$ have to be regarded as known functions. Thus,
with reference to (\ref{Hassnv})-(\ref{TNN}), the corresponding semi-discrete Hamiltonian problem is given by:
\begin{eqnarray}\label{temp}
\dbq &=& \bp ~\equiv~\frac{1}\Dx\nabla_\bp H,  \qquad t>0,\\ \nonumber
\dbp &=& -\frac{1}{\Dx ^2} T_N \bq+\frac{\Sigma }{\Dx^2}\bw(\bq,t)- f'(\bq ) ~\equiv~-\frac{1}\Dx\nabla_\bq H,
\end{eqnarray}
where
$$\Sigma = \pmatrix{ccccc} -1 \\ &0\\ &&\ddots \\ &&&0\\ &&&&1\endpmatrix\in\RR^{N\times N}.$$
By considering (\ref{u0Npiu1}), one has then
$$%\begin{equation}\label{bfiN}
\frac{\Sigma}{\Dx^2}\bw(\bq,t) = 
\frac{1}{\Dx^2} \pmatrix{c} 
u_0-u_1\\
0\\
\vdots
\\
0\\
u_{N+1}-u_N
\endpmatrix = \frac{1}{\Dx}
\pmatrix{c}-\phi_0(t)\\ 0\\ \vdots \\0\\ \phi_1(t)\endpmatrix \equiv \frac{1}{\Dx} \bffi(t),
$$%\end{equation}
thus obtaining the final shape of (\ref{temp}):
\begin{eqnarray}\label{wave1na}
\dbq &=& \bp,  \qquad t>0,\\
\dbp &=& -\frac{1}{\Dx ^2} T_N \bq+\frac{1}{\Dx}\bffi(t)- f'(\bq ).\nonumber
\end{eqnarray}
As in the case of Dirichlet boundary conditions, problem (\ref{wave1na}) is Hamiltonian with the non-autonomous Hamiltonian (\ref{Hassnv}): again, we can transform this latter into an autonomous one, by introducing the couple of auxiliary conjugate variables (\ref{auxvar}) and the augmented Hamiltonian (compare with (\ref{Hmod}))
\begin{equation}\label{Hmod1}
\tH(\bq,\bp,\tq,\tp) = H(\bq,\bp,\tq)+\tp,
\end{equation}
with $H$ now given by (\ref{Hassnv}).
The dynamical system corresponding to this new Hamiltonian function is, for $t>0$:
\begin{eqnarray}\nonumber
\dbq &=& \bp ~\equiv~\frac{1}\Dx\nabla_\bp \tH,   \\ \nonumber
\dbp &=& -\frac{1}{\Dx ^2} T_N \bq+\frac{1}{\Dx}\bffi- f'(\bq ) ~\equiv~-\frac{1}\Dx\nabla_\bq \tH,\\ \label{wave2n}
\dtq  &=& 1 ~\equiv~\frac{\partial}{\partial\tp} \tH,\\ \nonumber
\dtp  &=& -\frac{\partial}{\partial\tq} \tH,
\end{eqnarray}
with initial conditions as in (\ref{init1}). Concerning the last equation in (\ref{wave2n}), from (\ref{Hmod1}), (\ref{Hassnv})-(\ref{TNN}), and (\ref{u0Npiu1}), one has,
by considering that $q_i(t) \equiv u_i(t)$, $q_i'(t) = p_i(t)\equiv v_i(t)$ (see (\ref{xqp})), and $\tq\equiv t$,
\begin{eqnarray}\nonumber
\dtp  &=& -\frac{\partial}{\partial\tq} \tH ~=~ -\frac{1}{\Dx}\left( (u_0(\tq)-u_1)v_0(\tq) + (u_{N+1}(\tq)-u_N)v_{N+1}(\tq)\right)\\
&=& \phi_0(\tq)[v_1-\Dx\phi_0'(\tq)] - \phi_1(\tq)[v_N+\Dx \phi_1'(\tq)].\label{lasteq}
\end{eqnarray}
Now, problem (\ref{wave2n})-(\ref{lasteq}) is Hamiltonian with an autonomous Hamiltonian function, so that its energy  (\ref{Hmod1}) is conserved (clearly, considerations similar to those reported in Remark~\ref{prescri} for the Dirichlet case can be now repeated). 
 
 Also now, the discrete problem  can be cast in vector form, formally as done in (\ref{by1})-(\ref{wave2m}).
Moreover, concerning the discretization issue, arguments similar to those seen in Section~\ref{disass} apply to the present case. In particular, the following result holds true, the proof being similar to that of Theorems~\ref{hbvmpbc} and \ref{hbvmpbc1}.

\begin{theo}\label{hbvmpbc2} Let ~$\by_1=\bu(h)$~ be the approximation to $\by(h)$, solution of (\ref{by1})-(\ref{wave2m}),  with $\tH$ as in (\ref{Hmod1})-(\ref{lasteq}), provided by a HBVM$(k,s)$ method used with stepsize $h$. One then obtains:
$$\by_1-\by(h) = O(h^{2s+1}),$$ that is the method has order $2s$. Moreover, assuming that $f$, $\phi_0$, and $\phi_1$ in (\ref{wave})-(\ref{assnbc}) are suitably regular:
$$\tH(\by_1)-\tH(\by_0) = \left\{\begin{array}{ccl}
0, &\quad&\mbox{if ~$f\in\Pi_\nu$, $\phi_0,\phi_1\in\Pi_\rho,$~ with}\\[.2cm]
    &         & 2k\ge\max\{\nu s, 2\rho+s-1, 2s+\rho\},\\[.5cm]
O(h^{2k+1}), &&\mbox{otherwise.}
\end{array}\right.$$
\end{theo}
Evidently,  considerations similar to those stated in Remark~\ref{rem1} can be repeated also in the present situation.

\section{Periodic boundary conditions revisited}\label{pbcr}

The case of periodic boundary conditions, i.e. (\ref{wave})-(\ref{perbc}), deserves to be further investigated. In fact, the finite-difference discretizations considered above, turn out to provide a second-order spatial accuracy, in the used stepsize $\Dx$. When either Dirichlet or Neumann boundary conditions are specified, it is not possible to easily derive higher-order semi-discrete Hamiltonian formulations of the problem. Conversely, in the case of periodic boundary conditions, this can be easily accomplished.  As matter of fact, by suitably replacing the circulant matrix $T_N$ defined in (\ref{TNc}), one obtains that the Hamiltonian (\ref{Hperd}) remains formally the same, as well as the semi-discrete Hamiltonian problem (\ref{wave1d}). For this purpose, any symmetric high-order approximation to the second spatial derivative could be used (see e.g., \cite{AmSg05}), to derive a new circulant and symmetric band-matrix. As an example, the following matrix provides a fourth-order spatial approximation   \cite{QinZhang90},%\marginpar{\red $-T_N$}
\begin{equation}\label{TNc4}
T_N = -\pmatrix{rrrrrrrr}
-\frac{5}2     &\frac{4}3        &-\frac{1}{12}  &                   &            &-\frac{1}{12} &\frac{4}3\\
\frac{4}3      &\ddots            &\ddots           &\ddots          &            &               &-\frac{1}{12}\\
-\frac{1}{12} &\ddots           &\ddots           &\ddots           &\ddots\\
                    &\ddots          &\ddots           &\ddots            &\ddots &\ddots\\
                    &                    &\ddots           &\ddots             &\ddots &\ddots &-\frac{1}{12}\\
-\frac{1}{12}&                    &                     &\ddots             &\ddots &\ddots &\frac{4}3\\[1mm]
\frac{4}3    &-\frac{1}{12} &                     &                       &-\frac{1}{12}&\frac{4}3 &-\frac{5}2 
\endpmatrix\in\RR^{N\times N},
\end{equation}
whereas, the following one provides a sixth-order spatial approximation (see \cite{AmSg05} for additional examples):%\marginpar{\red $-T_N$}
\begin{equation}\label{TNc6}
T_N = -\pmatrix{rrrrrrrrr}
-\frac{49}{18}     &\frac{3}2        &-\frac{3}{20}  &\frac{1}{90}                   &            &\frac{1}{90}  &-\frac{3}{20} &\frac{3}2 \\
\frac{3}2      &\ddots            &\ddots           &\ddots          &  \ddots          &          &\frac{1}{90}          &-\frac{3}{20}\\
-\frac{3}{20} &\ddots           &\ddots           &\ddots           &\ddots           & \ddots          &                           &\frac{1}{90}\\
\frac{1}{90}                    &\ddots          &\ddots           &\ddots            &\ddots &\ddots &\ddots\\
                    & \ddots                   &\ddots           &\ddots             &\ddots &\ddots &  \ddots  & \frac{1}{90}\\
\frac{1}{90}&                    &\ddots                     &\ddots             &\ddots &\ddots & \ddots &-\frac{3}{20}\\
-\frac{3}{20}    &\frac{1}{90} &                     &\ddots                       &\ddots &\ddots  & \ddots &\frac{3}2\\[1mm]
\frac{3}2 &-\frac{3}{20}    &\frac{1}{90} &                     &\frac{1}{90} &-\frac{3}{20} &\frac{3}2 & -\frac{49}{18} 
\endpmatrix\in\RR^{N\times N},
\end{equation}

\subsection{Fourier space discretization}

An alternative approach, which we shall investigate in the sequel, is that of using a Fourier approximation  in space (see, e.g., \cite{Faou12}). For this purpose, let us consider the following complete set of orthonormal functions in $[0,1]$:
\begin{equation}\label{cnsn}
c_0(x)\equiv 1, \qquad c_k(x) = \sqrt{2}\cos(2k\pi x), \quad s_k(x) =\sqrt{2}\sin(2k\pi x), \qquad k=1,2,\dots,
\end{equation}
so that
\begin{equation}\label{orto}
\int_0^1 c_i(x)c_j(x)\dd x=\int_0^1 s_i(x)s_j(x)\dd x = \delta_{ij}, \qquad \int_0^1 c_i(x)s_j(x)\dd x=0, \qquad
\forall i,j.
\end{equation}
The following expansion of the solution of (\ref{wave})-(\ref{perbc}) is a slightly different way of writing the usual Fourier expansion in space:
\begin{eqnarray}\nonumber
u(x,t) &=& c_0(x)\gamma_0(t) +\sum_{n\ge 1} \left[c_n(x)\gamma_n(t)+s_n(x)\eta_n(t)\right] \\
\label{expu}
&\equiv&\gamma_0(t) +\sum_{n\ge 1}\left[ c_n(x)\gamma_n(t)+s_n(x)\eta_n(t)\right], \qquad x\in[0,1],\quad t\ge0,
\end{eqnarray}
with
$$%\begin{equation}\label{gammaeta}
\gamma_n(t) = \int_0^1 c_n(x)u(x,t)\dd x, \qquad \eta_n(t) = \int_0^1 s_n(x)u(x,t)\dd x,
$$%\end{equation}
which is allowed because of the periodic boundary conditions (\ref{perbc}). Consequently, by taking into account (\ref{orto}),  the first equation in (\ref{wave}) can be rewritten as: 
\begin{eqnarray}\nonumber
\ddgam_n(t) &=& -(2\pi n)^2 \gamma_n(t) \\ \nonumber
&&- \int_0^1 c_n(x)f'\left( \gamma_0(t) +\sum_{n\ge 1} \left[c_n(x)\gamma_n(t)+s_n(x)\eta_n(t)\right] \right)\dd x, \quad n\ge0,\\ \label{fourier1} \\ \nonumber
\ddeta_n(t) &=& -(2\pi n)^2 \eta_n(t)\\ \nonumber
&& - \int_0^1 s_n(x)f'\left( \gamma_0(t) +\sum_{n\ge 1}\left[ c_n(x)\gamma_n(t)+s_n(x)\eta_n(t)\right]  \right)\dd x, \quad n\ge1,
\end{eqnarray}
where the double dot denotes, as usual, the second time derivative. 
The initial conditions are clearly given by (see (\ref{wave})):
\begin{eqnarray}\nonumber
\gamma_n(0) = \int_0^1 c_n(x) \psi_0(x)\dd x, &\qquad& \eta_n(0) = \int_0^1 s_n(x) \psi_0(x)\dd x,\\[-1.5mm] \label{gameta0}\\[-1.5mm]
\dot\gamma_n(0) = \int_0^1 c_n(x) \psi_1(x)\dd x, &\qquad& \dot\eta_n(0) = \int_0^1 s_n(x) \psi_1(x)\dd x.\nonumber
\end{eqnarray}
By introducing the infinite vectors
\begin{eqnarray}\nonumber
\bom(x) &=& \pmatrix{cccccc} c_0(x)& c_1(x)& s_1(x) & c_2(x)& s_2(x)& \dots\endpmatrix^\top,\\
\label{qp} \\ \nonumber%%\qquad \bp(t)=\dot\bq(t),
\bq(t) &=& \pmatrix{cccccc} \gamma_0(t) & \gamma_1(t)& \eta_1(t) & \gamma_2(t)& \eta_2(t)& \dots\endpmatrix^\top,
\end{eqnarray}
the infinite matrix
\begin{equation}\label{D}
D = \pmatrix{cccccc} 0 \\ &(2\pi)^2 \\ &&(2\pi)^2\\ &&&(4\pi)^2\\ &&&&(4\pi)^2\\ &&&&& \ddots\endpmatrix,
\end{equation} and considering that (see (\ref{expu}))
\begin{equation}\label{expu1}
u(x,t) = \bom(x)^\top \bq(t),
\end{equation}
problem (\ref{fourier1}) can be cast in vector form as:
\begin{eqnarray} \label{fourier2} 
\dot\bq(t) &=& \bp(t), \qquad t>0, \\ \nonumber
\dot\bp(t) &=& -D\bq(t) - \int_0^1 \bom(x) f'(\bom(x)^\top \bq(t))\dd x,
\end{eqnarray}
with the initial conditions (\ref{gameta0}) written, more compactly, as
\begin{equation}\label{qp0f}
\bq(0) = \int_0^1\bom(x)\psi_0(x)\dd x, \qquad \bp(0) = \int_0^1\bom(x)\psi_1(x)\dd x.
\end{equation}

The following result then holds true.
\medskip

\begin{theo}\label{thH1} Problem (\ref{fourier2}) is Hamiltonian, with Hamiltonian
\begin{equation}\label{H2}
H(\bq,\bp) = \frac{1}2\bp^\top \bp + \frac{1}2\bq^\top D\bq + \int_0^1 f(\bom(x)^\top \bq)\dd x.
\end{equation}
This latter is equivalent to the Hamiltonian (\ref{H1}), via the expansion (\ref{expu})-(\ref{expu1}).
\end{theo}
\proof The first statement is straightforward, by considering that
$$\nabla_\bq f(\bom(x)^\top \bq)) = f'(\bom(x)^\top \bq)\bom(x).$$
The second statement then easily follows, by taking into account (\ref{expu1}), from the fact that, see (\ref{vequt}), (\ref{orto}), (\ref{expu}),  and (\ref{qp}):
\begin{eqnarray*}
\int_0^1 v(x,t)^2 \dd x &=& \int_0^1 u_t(x,t)^2\dd x ~=~ \int_0^1 \left(\dgam_0(t)+ \sum_{n\ge1}\left[ \dgam_n(t) c_n(x) + \deta_n(t) s_n(x)\right]\right)^2\dd x\\ &=& \dgam_0(t)^2+\sum_{n\ge1}\left[ \dgam_n(t)^2 +\deta_n(t)^2\right] ~\equiv\bp(t)^\top \bp(t),
\end{eqnarray*}
and
\begin{eqnarray*}
\int_0^1 u_x(x,t)^2 \dd x  &=& \int_0^1  \left(\sum_{n\ge1} 2\pi n\left[\eta_n(t) c_n(x) - \gamma_n(t) s_n(x)\right]\right)^2\dd x\\
&=& \sum_{n\ge1} (2\pi n)^2\left[\eta_n(t)^2+ \gamma_n(t)^2\right] ~=~ \bq(t)^\top D\bq(t).
\end{eqnarray*}
\QED

\subsection{Truncated Fourier-Galerkin approximation}\label{trunc}
In order to obtain a practical computational procedure, we truncate the infinite expansion (\ref{expu}) to a finite sum:\begin{equation}\label{expuN}
u(x,t) ~\approx~ \gamma_0(t) +\sum_{n=1}^N \left[c_n(x)\gamma_n(t)+s_n(x)\eta_n(t)\right] ~\equiv~ u_N(x,t),
\end{equation}
which converges more than exponentially with $N$ to $u$, if this latter is an analytical function.\footnote{We refer, e.g., to \cite{CHQZ88}, for a corresponding comprehensive error analysis.}
In other words, we look for an approximation to $u(x,t)$ belonging to the functional subspace (see (\ref{cnsn}))
\begin{equation}\label{VN}
{\cal V}_N = \mathrm{span}\left\{c_0(x),c_1(x),\dots,c_N(x), s_1(x),\dots,s_N(x)\right\}.
\end{equation}
Clearly, such a truncated expansion will not satisfy problem (\ref{wave})-(\ref{perbc}). Nevertheless, in the spirit of Fourier-Galerkin methods \cite{B01}, by requiring that the residual
$$R(u_N) := (u_N)_{tt} - (u_N)_{xx} + f'(u_N)$$
be orthogonal to ${\cal V}_N$, one obtains the {\em weak formulation} of problem (\ref{wave})-(\ref{perbc}), consisting in the following set of $2N+1$ differential equations, 
\begin{eqnarray}\nonumber
\lefteqn{\ddgam_n(t) ~=~ -(2\pi n)^2 \gamma_n(t) }\\ \nonumber
&&- \int_0^1 c_n(x)f'\left( \gamma_0(t) +\sum_{n=1}^N \left[c_n(x)\gamma_n(t)+s_n(x)\eta_n(t)\right] \right)\dd x, \quad n=0,\dots,N,\\ \label{fourier1N} \\ \nonumber
\lefteqn{\ddeta_n(t) = -(2\pi n)^2 \eta_n(t)}\\ \nonumber
&& - \int_0^1 s_n(x)f'\left( \gamma_0(t) +\sum_{n=1}^N\left[ c_n(x)\gamma_n(t)+s_n(x)\eta_n(t)\right]  \right)\dd x, \quad n=1,\dots,N,
\end{eqnarray}
approximating the leading ones in (\ref{fourier1}). By defining the finite vectors  in $\RR^{2N+1}$ (compare with (\ref{qp})),
\begin{eqnarray}\nonumber
\bom_N(x) &=& \pmatrix{cccccccc} c_0(x)& c_1(x)& s_1(x) & c_2(x)& s_2(x)& \dots & c_N(x)& s_N(x)\endpmatrix^\top,\\ 
\label{qpN}\\[-3mm] \nonumber  %\qquad \bp_N(t)=\dot\bq_N(t), \qquad
\bq_N(t) &=& \pmatrix{cccccccc} \gamma_0(t)& \gamma_1(t)& \eta_1(t) & \gamma_2(t)& \eta_2(t)& \dots & \gamma_N(t)& \eta_N(t)\endpmatrix^\top,\end{eqnarray}
the matrix (compare with (\ref{D}))
\begin{equation}\label{DN}
D_N = \pmatrix{cccccccc} 0 \\ &(2\pi)^2 \\ &&(2\pi)^2\\ &&&(4\pi)^2\\ &&&&(4\pi)^2\\ &&&&&\ddots \\ &&&&&& (2N\pi)^2\\ &&&&&&& (2N\pi)^2\endpmatrix \in\RR^{(2N+1)\times (2N+1)},
\end{equation} and considering that (compare with (\ref{expuN}))
\begin{equation}\label{expu1N}
u_N(x,t) = \bom_N(x)^\top \bq_N(t),
\end{equation}
the equations (\ref{fourier1N}), which have to be satisfied by (\ref{expu1N}), can be cast in vector form as:
\begin{eqnarray} \label{fourier2N}
\dot\bq_N(t) &=& \bp_N(t), \qquad t>0, \\ \nonumber 
\dot\bp_N(t) &=& -D_N\bq_N(t) - \int_0^1 \bom_N(x) f'(\bom_N(x)^\top \bq_N(t))\dd x,
\end{eqnarray}
for a total of $4N+2$ differential equations. Clearly, from (\ref{gameta0}) one obtains that the initial conditions for (\ref{fourier2N}) are given by:
\begin{equation}\label{ge0}
\bq_N(0) = \int_0^1 \bom_N(x)\psi_0(x)\dd x, \qquad  \bp_N(0) = \int_0^1 \bom_N(x)\psi_1(x)\dd x.
\end{equation}
The following result then easily follows by means of arguments similar to those used to prove Theorem~\ref{thH1}.

\begin{theo}\label{thH2} Problem (\ref{fourier2N}) is Hamiltonian, with Hamiltonian
\begin{equation}\label{H2N}
H_N(\bq_N,\bp_N) = \frac{1}2\bp_N^\top \bp_N + \frac{1}2\bq_N^\top D_N\bq_N + \int_0^1 f(\bom_N(x)^\top \bq_N)\dd x.
\end{equation}
\end{theo}

We observe that (\ref{H2N}) is equivalent to a truncated Fourier expansion of the Hamiltonian (\ref{H1}) (see also (\ref{H2})).
Moreover, it is worth mentioning that  using the initial conditions (\ref{ge0}), in place of (\ref{qp0f}), results in an error $e_N$, in the initial data, given by
\begin{eqnarray}\nonumber
e_N^2 &=&  \int_0^1 \left(\psi_0(x) - \bom_N(x)^\top\bq_N(0)\right)^2\dd x +
\int_0^1 \left(\psi_1(x) - \bom_N(x)^\top\bp_N(0)\right)^2\dd x
\\ \nonumber
&=& \sum_{n>N} \left[\int_0^1c_n(x)\psi_0(x)\dd x\right]^2+ \left[\int_0^1s_n(x)\psi_0(x)\dd x\right]^2+\\
&& \sum_{n>N}\left[\int_0^1c_n(x)\psi_1(x)\dd x\right]^2 + \left[\int_0^1s_n(x)\psi_1(x)\dd x\right]^2. \label{error0}
\end{eqnarray}
However, it must be stressed that, unlike the finite-difference case, both $e_N$ and the approximation (\ref{H2N}) to the continuous Hamiltonian, converge more than exponentially in $N$ ($e_N$ to 0, and $H_N$ to $H$), provided that the involved functions are analytical.

\subsection{Full discretization}\label{intx}

Since problem (\ref{fourier2N}) is Hamiltonian, with an autonomous Hamiltonian, this latter is conserved along the solution. Consequently, energy conserving methods can be conveniently used for its solution. In particular, Theorem~\ref{hbvmpbc} continues formally to hold for HBVM$(k,s)$ methods. However, the integral appearing in (\ref{fourier2N}) need to be, in turn, approximated by means of a suitable quadrature rule. For this purpose, it could be convenient to do this by means of a composite trapezoidal rule, due to the fact that the argument is a periodic function. Consequently, having set
\begin{equation}\label{g_N}
\bg_N(x,t) = \bom_N(x) f'(\bom_N(x)^\top \bq_N(t)),
\end{equation} 
the uniform mesh on $[0,1]$
\begin{equation}\label{xi}
x_i = i\Delta x, \qquad i=0,\dots,m,\qquad \Delta x = \frac{1}m,
\end{equation}
and considering that 
~$\bg_N(0,t) = \bg_N(1,t),$~
one obtains:
\begin{eqnarray}\nonumber
\int_0^1 \bg_N(x,t)\dd x &=& \Delta x\sum_{i=1}^m \frac{\bg_N(x_{i-1},t)+\bg_N(x_i,t)}2 ~+~R(m)\\ 
&=& \frac{1}m\sum_{i=0}^{m-1} \bg_N(x_i,t) ~+~ R(m).\label{compo}
\end{eqnarray}
Let us study the error $R(m)$. For this purpose, we need some preliminary result.

\begin{lem}\label{lem1}
Let us consider the trigonometric polynomial 
\begin{equation}\label{poly}
p(x) = \sum_{k=0}^K \left[a_k\cos(2k\pi x)+ b_k \sin(2k\pi x)\right],
\end{equation}
and the uniform mesh (\ref{xi}). Then, for all $m\ge K+1$, one obtains: $$\int_0^1 p(x)\dd x = \frac{1}m\sum_{i=0}^{m-1} p(x_i).$$
\end{lem}
\proof See, e.g., \cite[Th.\,5.1.4]{DaBi08}.\,\QED
\bigskip

\begin{lem}\label{lem2}
Let us consider the trigonometric polynomial (\ref{poly}) and the uniform mesh (\ref{xi}).
Then,  for all $m\ge N+K+1$, one obtains: 
\begin{eqnarray}\label{cosint}
\int_0^1 \cos(2j\pi x) p(x)\dd x &=& \frac{1}m\sum_{i=0}^{m-1} \cos(2j\pi x_i) p(x_i),\\ \label{sinint}
\int_0^1 \sin(2j\pi x) p(x)\dd x &=& \frac{1}m\sum_{i=0}^{m-1} \sin(2j\pi x_i) p(x_i),\qquad j=0,\dots,N.
\end{eqnarray}
\end{lem}
\proof
By virtue of the prosthaphaeresis formulae, one has, for all $j=0,\dots,N$ and $k=0,\dots,K$:
\begin{eqnarray*}
\cos(2j\pi x)\cos(2 k\pi x) &=& \frac{1}2\left[ \cos(2(k+j)\pi x) + \cos(2(k-j)\pi x)\right],\\
\cos(2j\pi x)\sin(2 k\pi x)  &=&  \frac{1}2\left[ \sin(2(k+j)\pi x) + \sin(2(k-j)\pi x) \right],\\
\sin(2j\pi x)\cos(2 k\pi x) &=& \frac{1}2\left[ \sin(2(k+j)\pi x) - \sin(2(k-j)\pi x)\right],\\
\sin(2j\pi x)\sin(2 k\pi x)  &=&  \frac{1}2\left[ \cos(2(k-j)\pi x) - \cos(2(k+j)\pi x) \right].
\end{eqnarray*}
Consequently, the integrals at the left-hand side in (\ref{cosint})-(\ref{sinint}) are trigonometric polynomials of degree at most $N+K$. By virtue of Lemma~\ref{lem1}, it then follows that they are exactly computed by means of the composite trapezoidal rule at the corresponding right-hand sides, provided that $m\ge N+K+1$.\QED 
\bigskip

By virtue of Lemma~\ref{lem2}, the following result follows at once.

\begin{theo}\label{esatti}
Let the function $f$ appearing in (\ref{g_N}) (see also (\ref{expu1N})) be a polynomial of degree $\nu$, and let us consider the uniform mesh (\ref{xi}). Then, with reference to (\ref{compo}), for all $m\ge \nu N+1$ one obtains:
$$R(m)= 0\qquad i.e.,\qquad \int_0^1 \bg_N(x,t)\dd x = \frac{1}m\sum_{i=0}^{m-1} \bg_N(x_i,t).$$
\end{theo}

For a general function $f$, the following result holds true.

\begin{theo}\label{approssimati}
Let the function $\bg_N(x,t)$ defined at (\ref{g_N}), with $t$ a fixed parameter, belong to $W_{per}^{r,p}$, the Banach space of periodic functions on $\RR$ whose distribution derivatives up to order $r$ belong to $L_{per}^p(\RR)$. Then, with reference to (\ref{xi})-(\ref{compo}), one has:
$$ R(m)  = O(m^{-r}).$$
\end{theo}
\proof See \cite[Th.\,1.1]{KuRa09}.\QED\bigskip

We end this section by mentioning that different approaches could be also used for approximating the integral appearing in (\ref{fourier2N}):  we refer, e.g., to \cite{EvWe99}, for a comprehensive review on this topic.

\section{Implementation of the methods}\label{implement}
The efficient implementation of HBVMs has been studied in \cite{BIT11,BFI13,BFI13_2}. We here sketch the application of a HBVM$(k,s)$ method for solving (\ref{wave1d}), since the application to (\ref{wave1m}), (\ref{wave2n}), and (\ref{fourier2N}) is similar.
We consider the very first application of the method, so that the index of the time step can be skipped. As remarked in \cite{BIT11}, the discrete problem generated by a HBVM$(k,s)$ method is more conveniently recast in terms of the $s$  coefficients of the polynomial (\ref{dotu}), instead of the $k$ stages of the Runge-Kutta formulation (\ref{RKform}). Moreover, since in the case of the semi-discrete formulation of the wave equation the Hamiltonian is separable, additional savings are possible, since the dimension of the problem can be halved, as we are going to sketch.\footnote{This is not the case when considering different Hamiltonian PDEs, such as, e.g., the nonlinear Schr\"odinger equation.}
Let us then split the stage vector $Y$ of the Runge-Kutta formulation, into $Q$ and $P$, corresponding to the stages for $\bq$ and $\bp$, respectively. Consequently, from (\ref{RKform}) and (\ref{wave1d})-(\ref{q0p0}), one obtains, by setting $\bq_0=\bq(0)$,  $\bp_0=\bp(0)$, and $h$ the time step: 
\begin{equation}\label{QPform}
Q = \be\otimes \bq_0 + h\calI\calP^\top\Omega\otimes I_N\,P,\qquad P=\be\otimes \bp_0-h\calI\calP^\top\Omega \otimes I_N\, F(Q),
\end{equation}
where (see (\ref{wave1d}) and (\ref{TNc}))
\begin{equation}\label{FQ}
F(Q) = \frac{1}{\Dx^2}I_k\otimes T_N\, Q + f'(Q),
\end{equation}
with an obvious meaning of $f'(Q)$. By considering the following properties of the matrices $\calP$ and $\calI$, due to corresponding properties of Legendre polynomials \cite{BIT11},
\begin{itemize}
\item $\calI\calP^\top\Omega\be = \bc,$
\item
$\calP^\top\Omega\calI = X_s\equiv 
\pmatrix{cccc} 
\frac{1}2 &-\xi_1\\
\xi_1 &0 &\ddots\\
         &\ddots &\ddots &-\xi_{s-1}\\
         &            &\xi_{s-1}   &0\endpmatrix$, \\ with
         \begin{equation}\label{xi_i} \xi_i=\left(2\sqrt{4i^2-1}\right)^{-1},\qquad i=1,\dots,s-1,\end{equation}
\end{itemize}                 
%and ~$\calI\calP^\top\Omega\be = \bc,$~ 
substitution of the latter equation in (\ref{QPform}) in the former one gives:
$$Q = \be\otimes \bq_0 + h\bc\otimes \bp_0 -h^2\calI X_s\calP^\top\Omega\otimes I_N\, F(Q).$$          
By setting\,\footnote{Here, $\gamma_j$ is given by the entries of the vector $\hat\gamma_j(\bu)$ in (\ref{dotu}) corresponding to the $\bq$ components only. Consequently, it has a halved dimension, w.r.t. this latter vector.}
$$\bgam =  \calP^\top\Omega\otimes I_N\, F(Q) \equiv \pmatrix{c}
\gamma_0\\ \vdots\\ \gamma_{s-1}\endpmatrix,$$ 
one then obtains the following discrete problem (of block dimension $s$):
\begin{equation}\label{gamma}
G(\bgam) \equiv \bgam -\calP^\top\Omega\otimes I_N\, F\left(\be\otimes \bq_0 + h\bc\otimes \bp_0 -h^2\calI X_s\otimes I_N\,\bgam\right)=\bf0.\end{equation}           
Once (\ref{gamma}) is solved, the new approximations are then given by (see (\ref{xi_i}))  \cite{BIT11}:
$$\bp_1 = \bp_0 + h\gamma_0, \qquad \bq_1 = \bq_0 + h\bp_0 +h^2\left( \frac{1}2\gamma_0 -\xi_1\gamma_1\right).$$
Consequently, the solution of the discrete problem (\ref{gamma}) is the bulk of the computational cost of the step. For its solution, one could use the following simplified Newton iteration,
\begin{equation}\label{etaell}
\left(I_s\otimes I_N +\frac{h^2}{\Dx^2} X_s^2\otimes T_N\right)\Delta \bgam^\ell = -G(\bgam^\ell) \equiv \bet^\ell, \qquad \ell=0,1,\dots,
\end{equation}
which only considers the (main) linear part of the function $F$ (see (\ref{FQ})). However, even though the coefficient matrix of such iteration is constant, nevertheless, it has dimension $sN$. To reduce the computational cost, it is then better to use a {\em blended iteration} \cite{BIT11} (see also \cite{Br00,BrMa02,BrMa04}), formally defined as:
\begin{eqnarray}\label{blendit}
\bet_1^{\ell} &=& \rho_s^2 X_s^{-2}\otimes I_N\, \bet^{\ell},\\ \label{blendit1}
\Delta\bgam^\ell &=& I_s\otimes M_N^{-1}\left[ \bet_1^{\ell} + I_s\otimes M_N^{-1}\left( \bet^{\ell}-\bet_1^{\ell}\right)\right],
\qquad \ell=0,1,\dots,
\end{eqnarray}
where 
$$\rho_s = \min_{\lambda\in\sigma(X_s)} |\lambda|, \qquad M_N = I_N+\left( \frac{h\rho_s}{\Dx}\right)^2T_N,$$
with $\sigma(X_s)$ denoting the spectrum of matrix $X_s$. Consequently, the computational cost of each iteration is given by:
\begin{description}
\item{- the evaluation of $\bet^\ell$ in (\ref{etaell}).} This requires $k$ evaluations of the right-hand side of the second equation in (\ref{wave1d}) (see (\ref{FQ})--(\ref{etaell})) plus
$(4ks+3k+s)N$ {\em flops};\footnote{We count as 1 {\em flop}, one elementary {\em fl}oating-point {\em op}eration.}

\item{- the evaluation of $\bet_1^\ell$ in (\ref{blendit}).}  Concerning matrix $\rho_s^{-1}X_s$, one can either invert and square it in advance, so that the costs for computing $\bet_1^\ell$ is $2s^2N$ {\em flops},
or solve 2 tridiagonal linear systems, so that, once the factorization is computed,\footnote{This costs less than $3s$ {\em flops}.} the cost per iteration amounts to $10sN$ {\em flops}.
Consequently, the corresponding computational cost is given by $2\min\{s,5\}sN$ {\em flops};

\item{- the evaluation of $\Delta\bgam^\ell$ in (\ref{blendit1}).}  This requires solution of $2s$ linear systems with the symmetric matrix $M_N$ plus $2sN$ {\em flops}. Concerning  matrix $M_N$, an additional saving of computational effort is gained by retaining only its tridiagonal part (or by considering an approximate inverse).\footnote{In general, the matrix becomes {\em banded}, when considering higher-order discretizations, see, e.g., (\ref{TNc4})-(\ref{TNc6}).} In such a case, after its factorization,\footnote{This costs less than $3N$ {\em flops}.} one has a cost of less than 
$10sN$ {\em flops}. The total cost is then less than $12sN$ {\em flops}.
\end{description}
 In conclusion, the total cost per iteration amounts to 
 $k$ function evaluations plus 
$(13s+3k+2\min\{s,5\}s+4ks)N$ {\em flops}.

It is worth mentioning that the same complexity is obtained in the case of Dirichlet or Neumann boundary conditions, by considering the corresponding tridiagonal matrices (\ref{TND}) and (\ref{TNN}), respectively. Instead, when using the Fourier-Galerkin spatial semi-discretization, one obtains that matrix $M_N$ is given by
$$M_N = I_{2N+1} +(h\rho_s)^2 D_N \in\RR^{(2N+1)\times (2N+1)},$$
where matrix $D_N$ is {\em diagonal} (see (\ref{DN})). Consequently, also $M_N$ is a diagonal matrix and, therefore, the complexity per iteration, besides the functions evaluations of the second equation in (\ref{fourier2N}) (which are the same as before i.e., $k$), decreases. As matter of fact, the required {\em flops} per iteration are now given by the dimension of the problem, times a factor $(5s+3k+2\min\{s,5\}s+4ks)$, in place of the factor $(13s+3k+2\min\{s,5\}s+4ks)$ seen above.

As a result of the previous arguments, one then expects a complexity per step which is {\em linear} in the dimension of the problem and, therefore, comparable with that of an explicit method. Moreover, in contrast to the $A$-stable HBVM$(k,s)$ methods,
explicit methods may suffer from stepsize restrictions due to stability reasons, as we shall see in the numerical tests.

\section{Numerical tests}\label{numtest}

\begin{figure}
\centerline{\includegraphics[width=9cm,height=7.5cm]{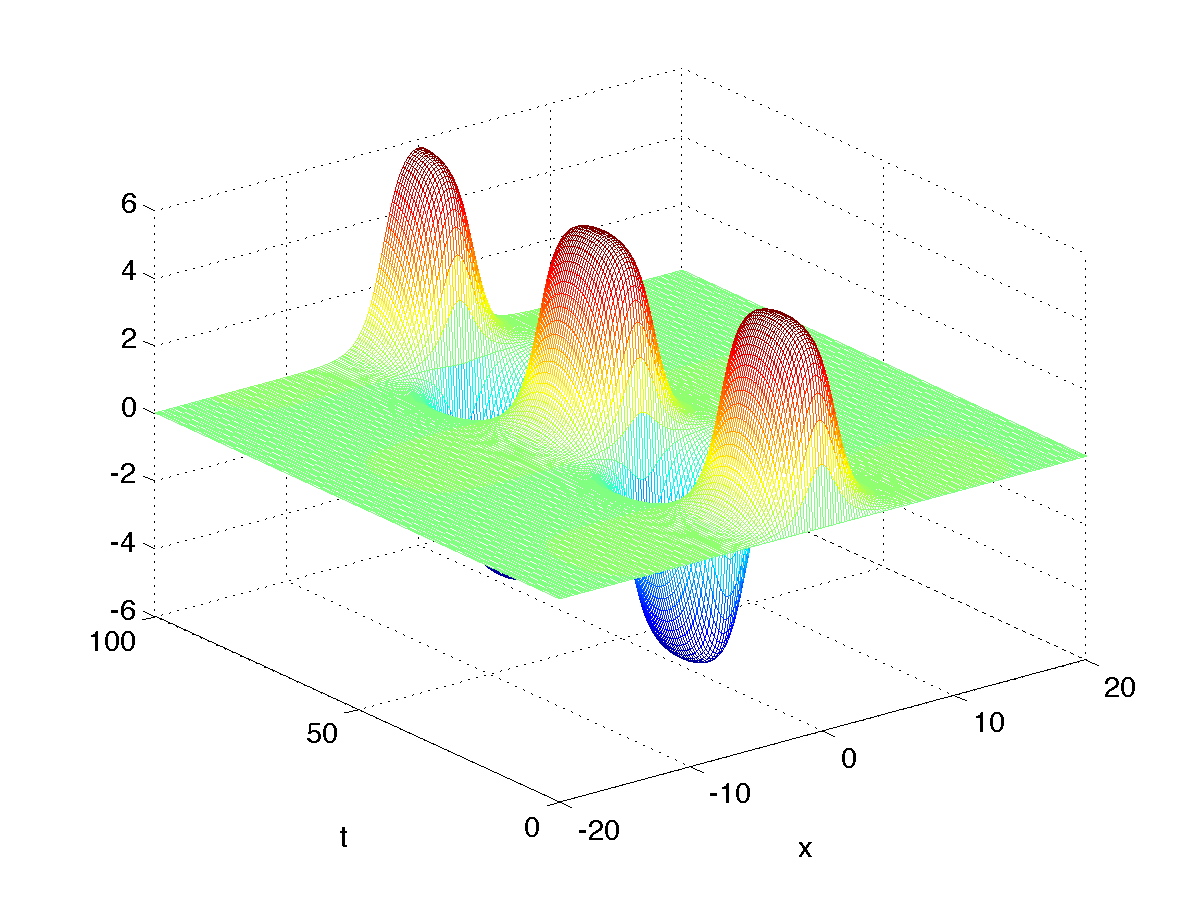}\quad\includegraphics[width=9cm,height=7.5cm]{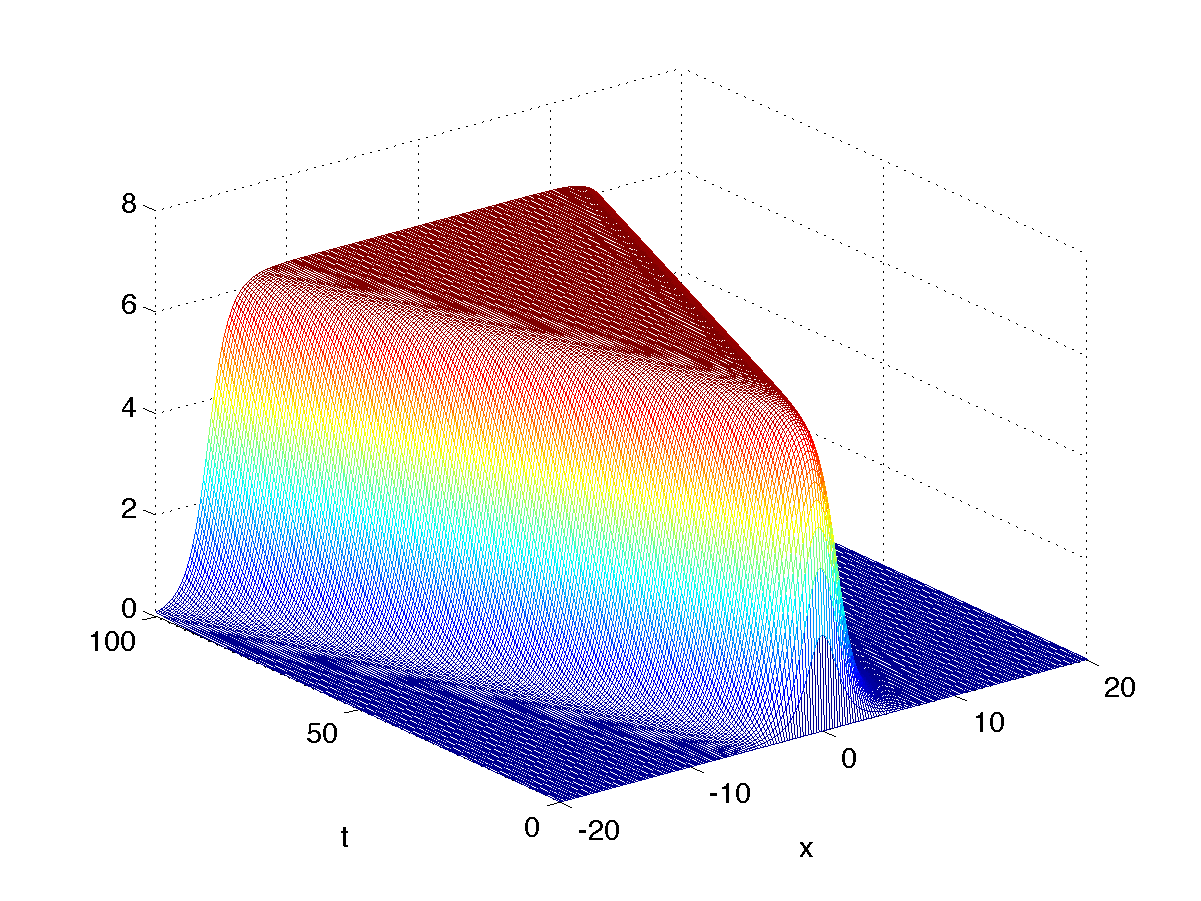}}
\caption{{\em Breather} solution for $\gamma=1.01$ (left plot), and {\em kink-antikink} solution for $\gamma=0.99$ (right plot).}
\label{soli0}
\end{figure}

\begin{figure}
\centerline{\includegraphics[width=9cm,height=7.5cm]{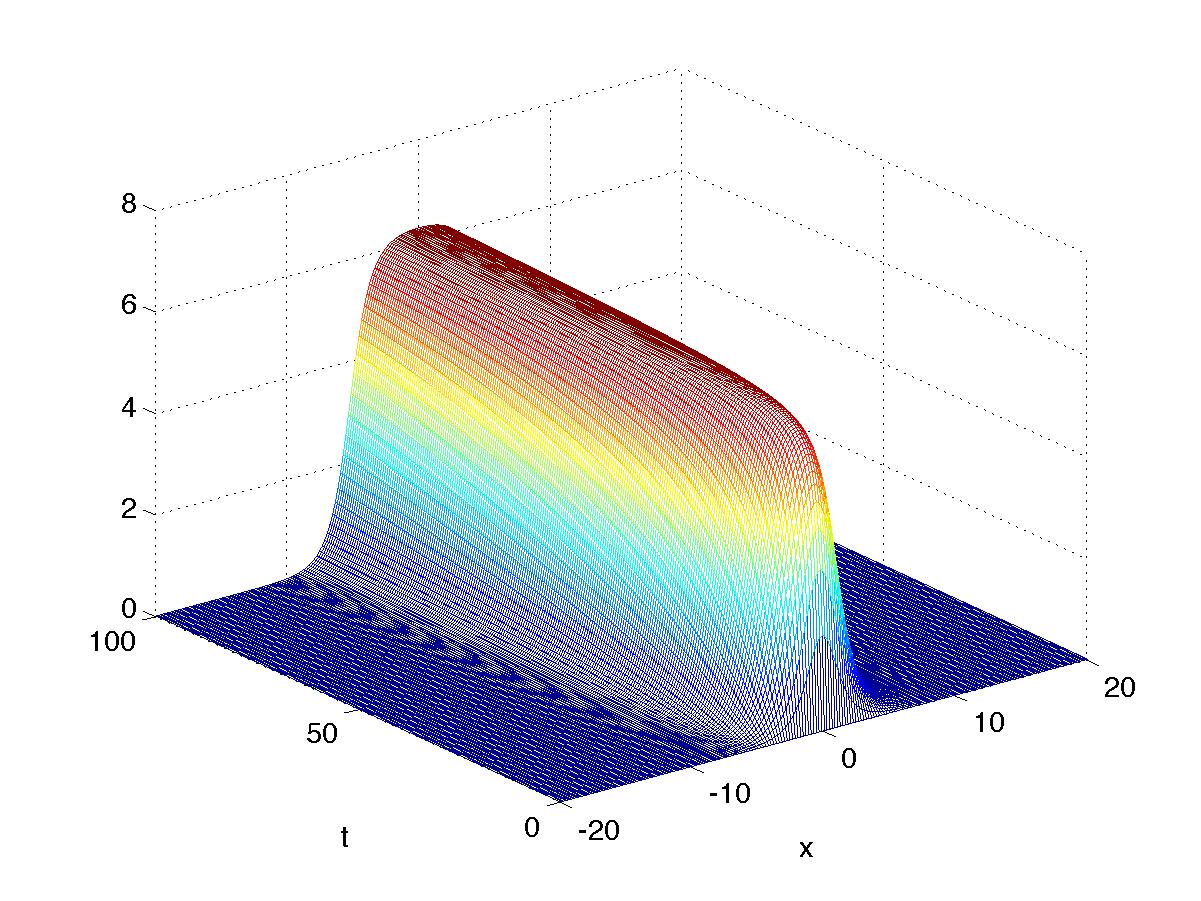}\quad\includegraphics[width=9cm,height=7.5cm]{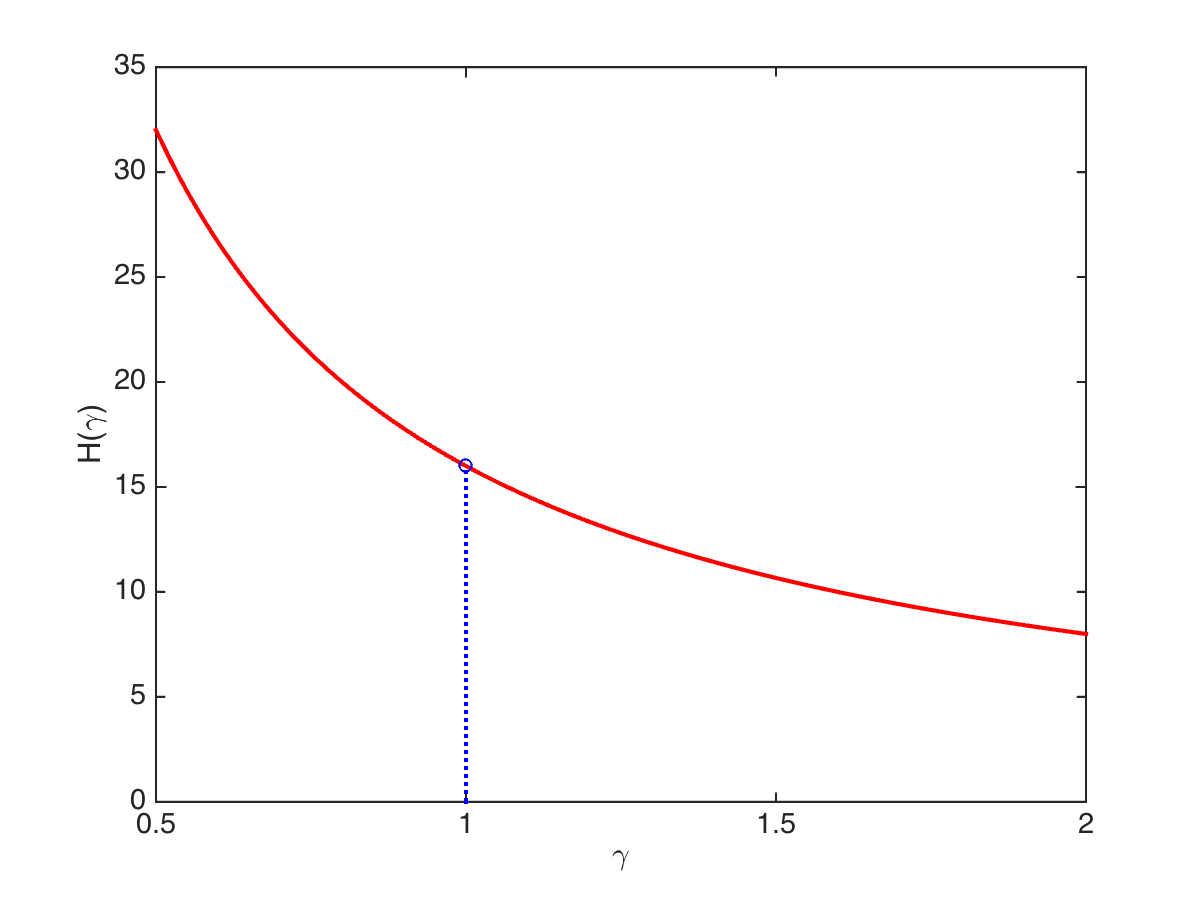}}
\caption{{\em Double-pole} solution for $\gamma=1$ (left plot), and Hamiltonian as a function of $\gamma$ (right plot).}
\label{soli1}
\end{figure}

We here consider a few numerical tests, concerning the so called {\em sine-Gordon} equation, which is in the form (\ref{wave}):
\begin{equation}\label{sineG}
u_{tt}(x,t) =u_{xx}(x,t)-\sin(u(x,t)), \qquad x\in[-20,20], \quad t\ge0.
\end{equation}
In particular, we shall consider {\em soliton-like} solutions, as described in \cite{W07}, defined by the initial conditions:
\begin{equation}\label{sineG0}
u(x,0) \equiv 0, \qquad u_t(x,0) = \frac{4}\gamma \sech\left(\frac{x}\gamma\right), \qquad \gamma>0.
\end{equation}
Depending on the value of the positive parameter $\gamma$, the solution is known to be given by:
\begin{equation}\label{sineGu}
u(x,t) = 4\atan\left[ \varphi(t;\gamma)\, \sech\left(\frac{x}\gamma\right)\right], 
\end{equation}
with
\begin{equation}\label{fi}
\varphi(t;\gamma) = \left\{ \begin{array}{lcc}
(\sqrt{\gamma^2-1})^{-1} \sin\left( \gamma^{-1}\sqrt{\gamma^2-1} t\right), &\mbox{if}~ &\gamma>1,\\[2mm]
 ~t, &\mbox{if}~  &\gamma=1,\\[2mm]
(\sqrt{1-\gamma^2})^{-1} \sinh\left( \gamma^{-1}\sqrt{1-\gamma^2} t\right), &\mbox{if}~  &0<\gamma<1.\\
\end{array}\right.
\end{equation}
The three cases are shown in Figures~\ref{soli0}--\ref{soli1}: on the left of Figure~\ref{soli0} is the plot of 
the first soliton (obtained for $\gamma>1$), which is named {\em breather}\/; on the right plot of Figure~\ref{soli0} is the case $0<\gamma<1$, which is named {\em kink-antikink}\/; at last, the case $\gamma=1$, which  is named {\em double-pole}, separates the two different types of dynamics and is shown in the left plot of Figure~\ref{soli1}.
Moreover, the space interval being fixed,\footnote{I.e., $[-20,20]$, in our case (see (\ref{sineG})).} the Hamiltonian is a decreasing function of $\gamma$, as is shown in the right plot of Figure~\ref{soli1}. This means that the value of the Hamiltonian characterizes the dynamics. Consequently, in a neighbourhood of $\gamma=1$, where the Hamiltonian assumes a value $\simeq 16$, nearby values of the Hamiltonian will provide different types of soliton solutions. As a result, energy conserving methods are expected to be useful, when numerically solving problem (\ref{sineG})-(\ref{sineG0}) with $\gamma=1$.

Let us then solve such a problem, at first with periodic boundary conditions, by using:
\begin{itemize}
\item a finite-difference approximation with $N=400$ equispaced mesh points;

\item a trigonometric polynomial approximation of degree $N=100$ and, moreover, $m=200$ equispaced mesh points.\footnote{In fact, $m=200$ is an appropriate choice for $N=100$, in this case.}
In so doing, the error (\ref{error0})  in the initial condition is $e_N\simeq 1.6\cdot 10^{-11}$, so that it is quite well matched.
\end{itemize}
For the time integration, let us consider the following second-order methods, used with stepsize $h=10^{-1}$ for $10^3$ integration steps:
\begin{itemize}
\item the (symplectic) implicit mid-point rule, i.e., HBVM(1,1), for which the Hamiltonian error is $\simeq 2\cdot 10^{-2}$ (though without a drift);

\item the (practically) energy-conserving HBVM(5,1) method, for which the Hamiltonian error is $\simeq 9\cdot 10^{-14}$.
\end{itemize}
 
\begin{figure}
\centerline{\includegraphics[width=9cm,height=7.5cm]{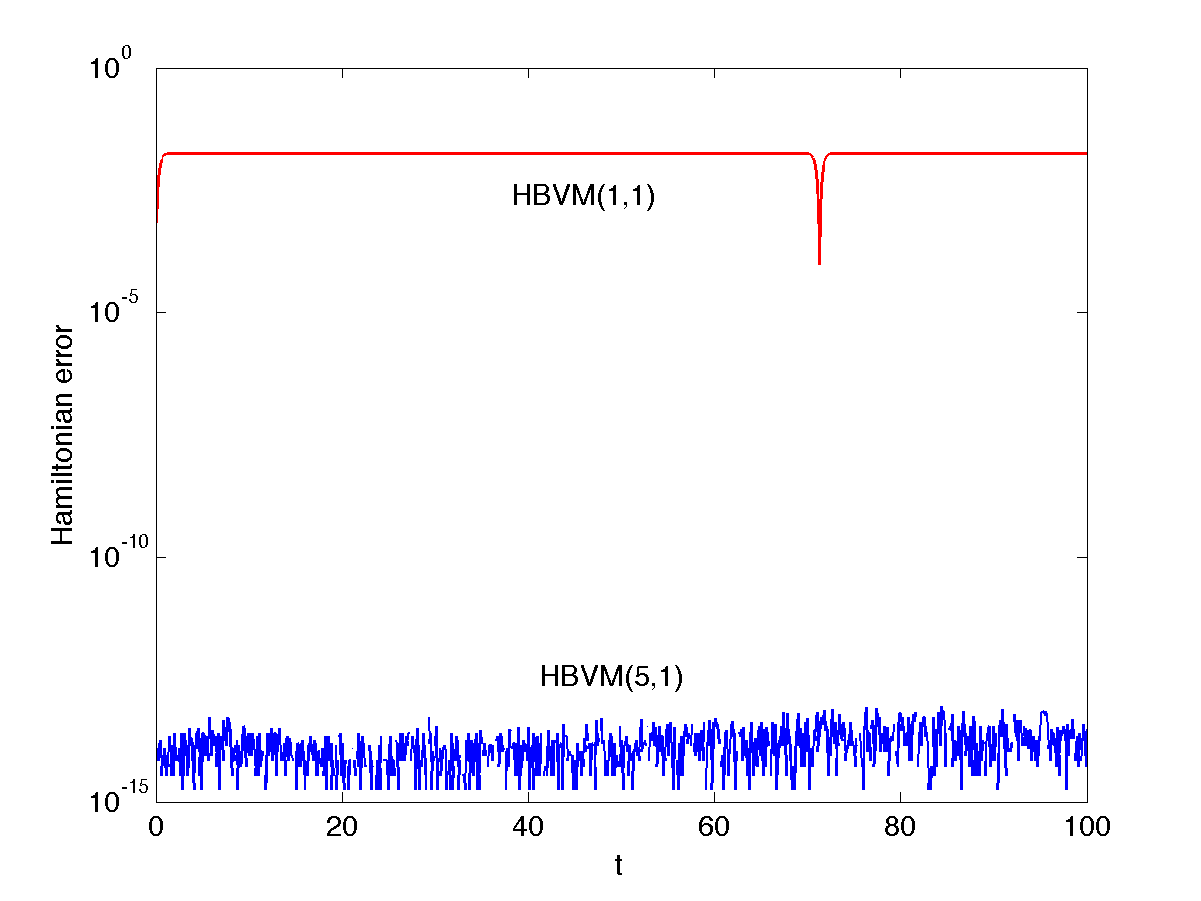}\quad\includegraphics[width=9cm,height=7.5cm]{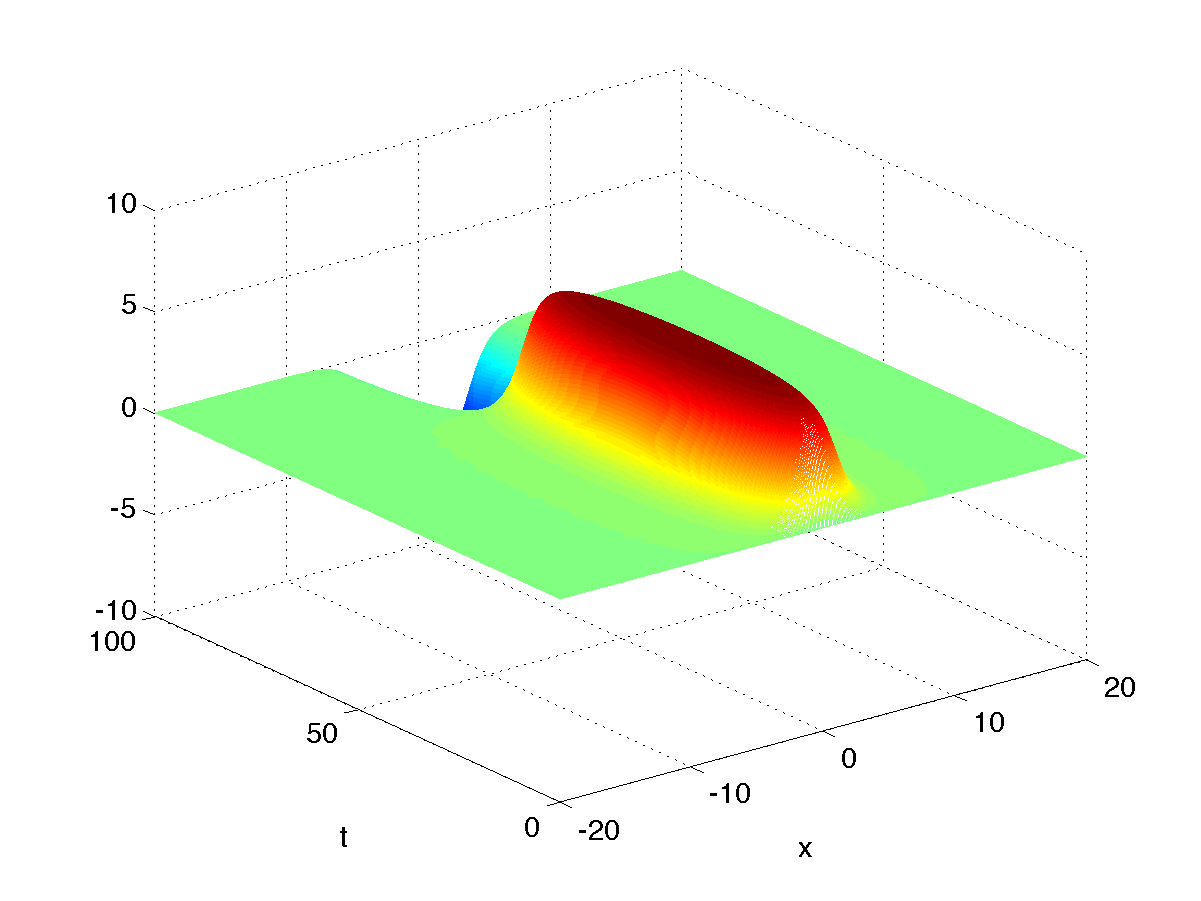}}
\caption{periodic boundary conditions and finite-difference approximation. Hamiltonian error (left plot) when using the HBVM(1,1) and HBVM(5,1) methods with stepsize $h=0.1$, and numerical solution provided by HBVM(1,1) (right plot) when solving problem (\ref{sineG})-(\ref{sineG0}) with $\gamma=1$.}
\label{test1}
\end{figure}

\begin{figure}[t]
\centerline{\includegraphics[width=9cm,height=7.5cm]{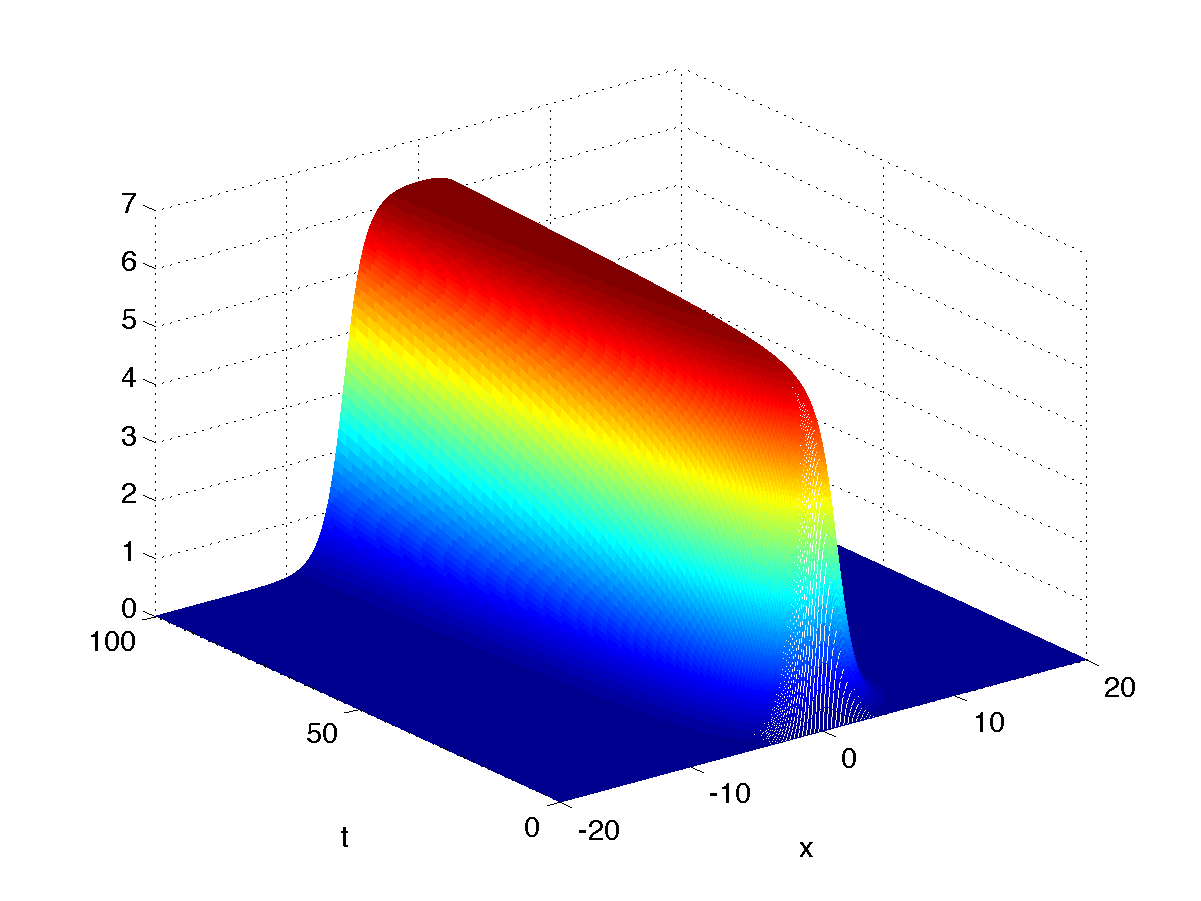}\quad\includegraphics[width=9cm,height=7.5cm]{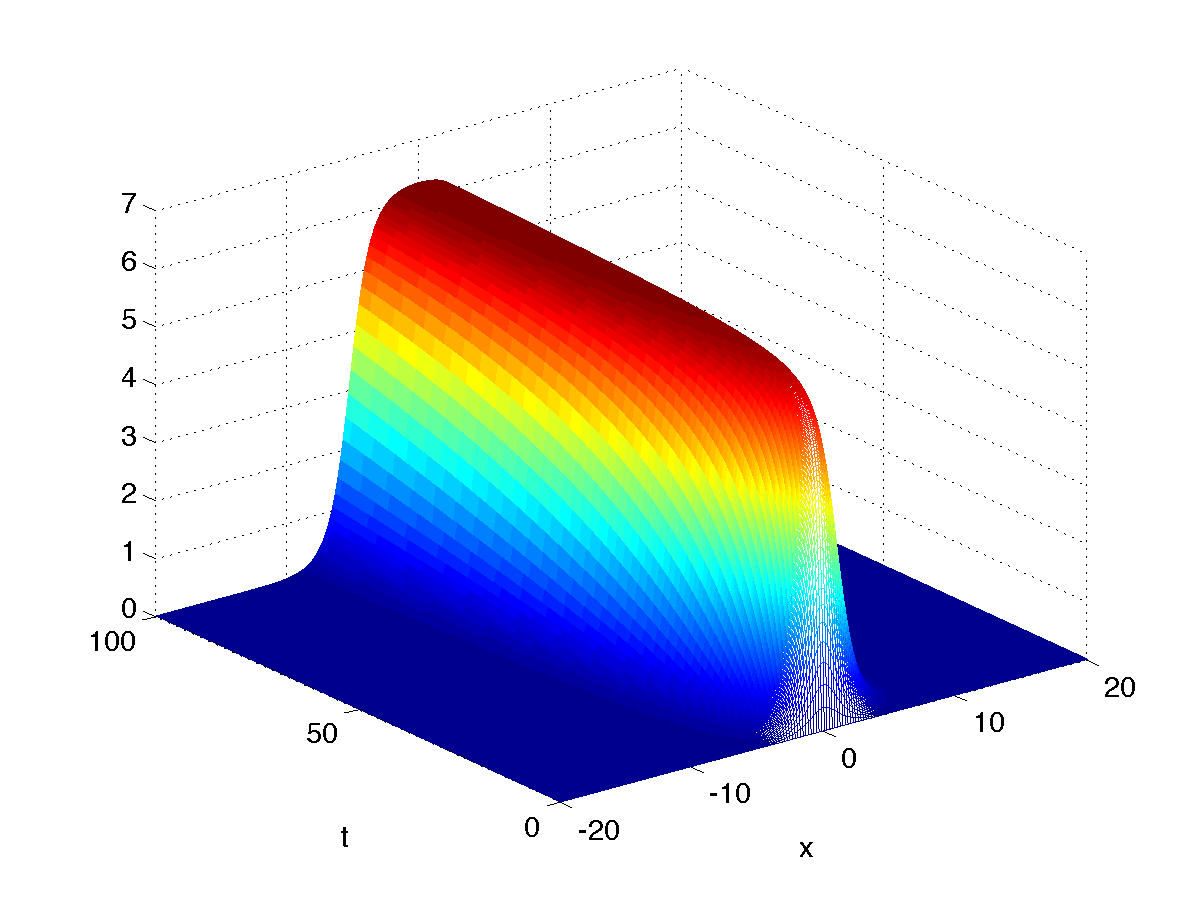}}
\caption{periodic boundary conditions. Computed solution by HBVM(5,1) with stepsize $h=0.1$ by using a finite-difference spatial discretization (left plot) or a spectral space discretization (right plot).}
\label{test1_1}
\end{figure}

Concerning the finite-difference space approximation, the error in the numerical Hamiltonian is plotted on the left of Figure~\ref{test1}.  The right plot of the same figure illustrates  the numerical approximation to the solution computed by the HBVM(1,1) method: as is clear, the computed approximation is wrong, since the method has provided a {\em breather}-like solution. On the contrary, HBVM(5,1) provides a correct approximation, qualitatively similar to that in the left-plot of Figure~\ref{soli1}: it is shown in the left plot in Figure~\ref{test1_1}.

\begin{figure}
\centerline{\includegraphics[width=9cm,height=7.5cm]{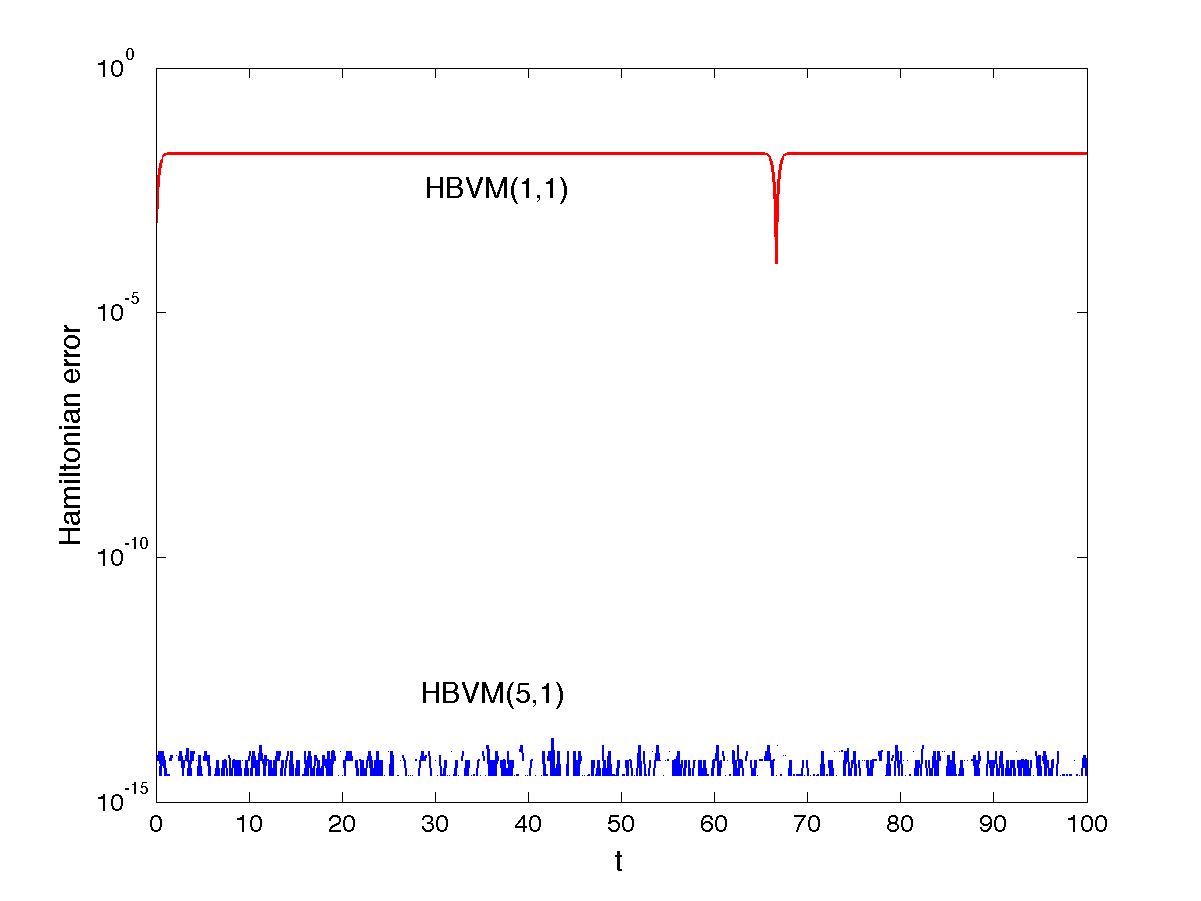}\quad\includegraphics[width=9cm,height=7.5cm]{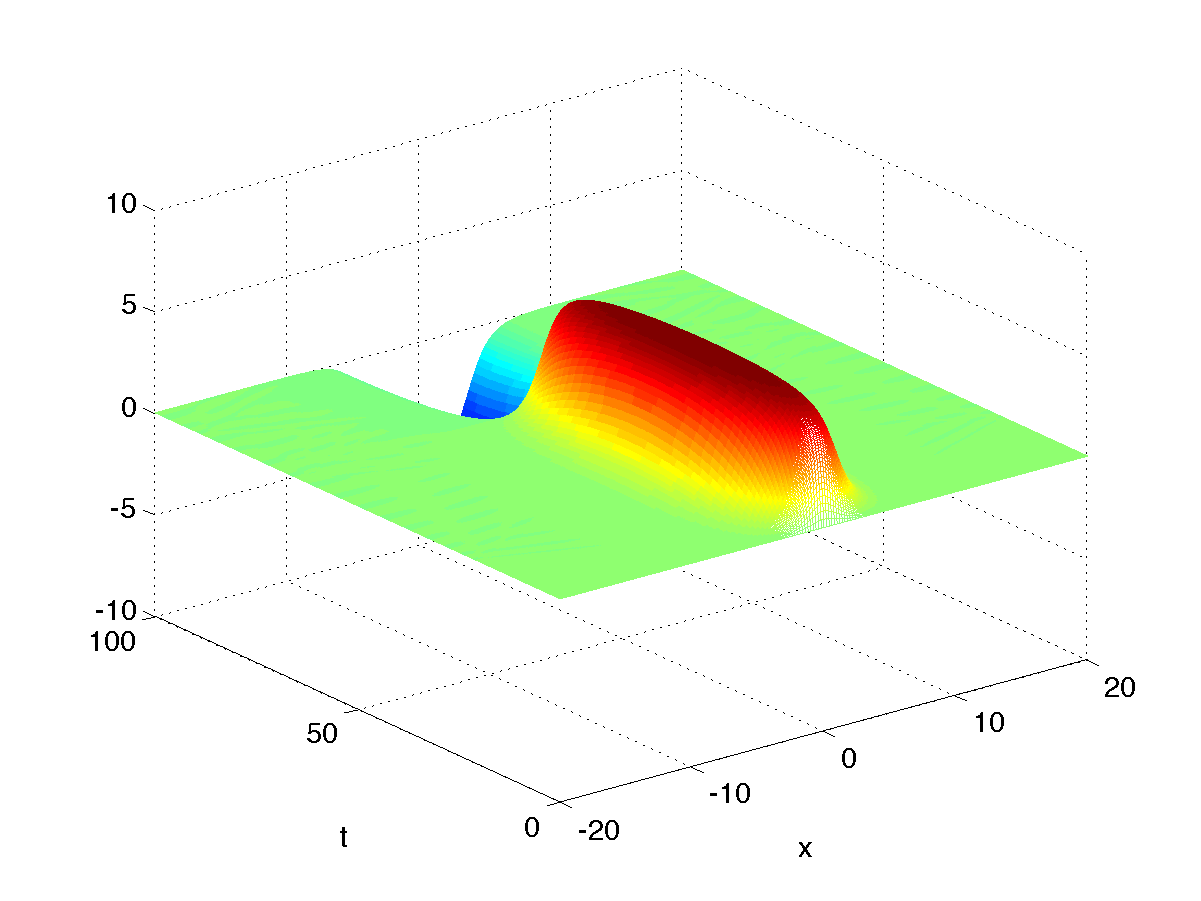}}
\caption{periodic boundary conditions and Fourier-Galerkin approximation. Hamiltonian error (left plot) when using the HBVM(1,1) and HBVM(5,1) methods with stepsize $h=0.1$, and numerical solution provided by HBVM(1,1) (right plot) when solving problem (\ref{sineG})-(\ref{sineG0}) with $\gamma=1$.}
\label{test1_2}
\end{figure}

Concerning the trigonometric polynomial approximation, the error in the numerical Hamiltonian is plotted on the left of Figure~\ref{test1_2}. The right plot of the same figure illustrates the numerical approximation to the solution computed by the HBVM(1,1) method:  it has again a {\em breather}-like shape and, thus, it is not qualitatively correct.  On the
contrary, HBVM(5,1) is able to reproduce the correct behaviour of the solution, as is
shown in the right plot in Figure~\ref{test1_1}.

\begin{figure}[t]
\centerline{\includegraphics[width=9cm,height=7.5cm]{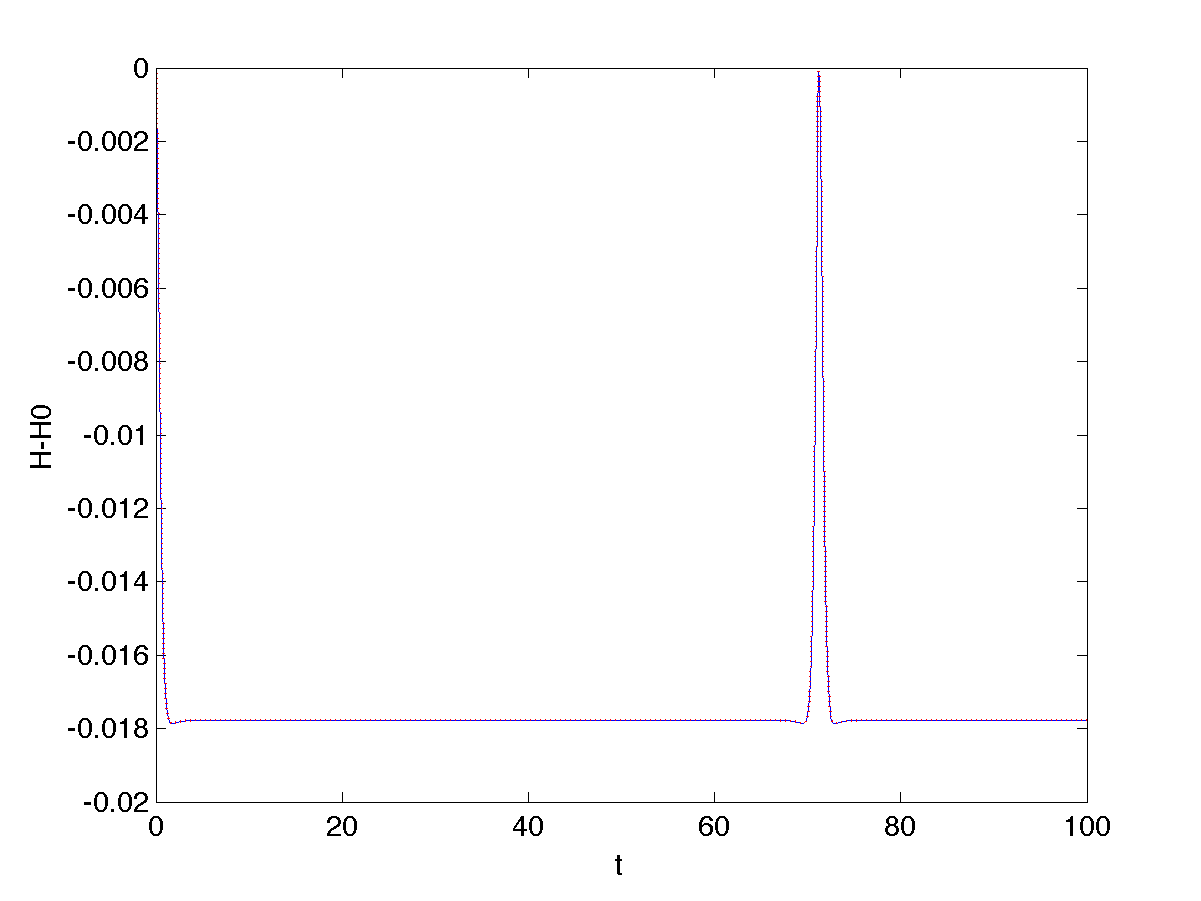}\quad\includegraphics[width=9cm,height=7.5cm]{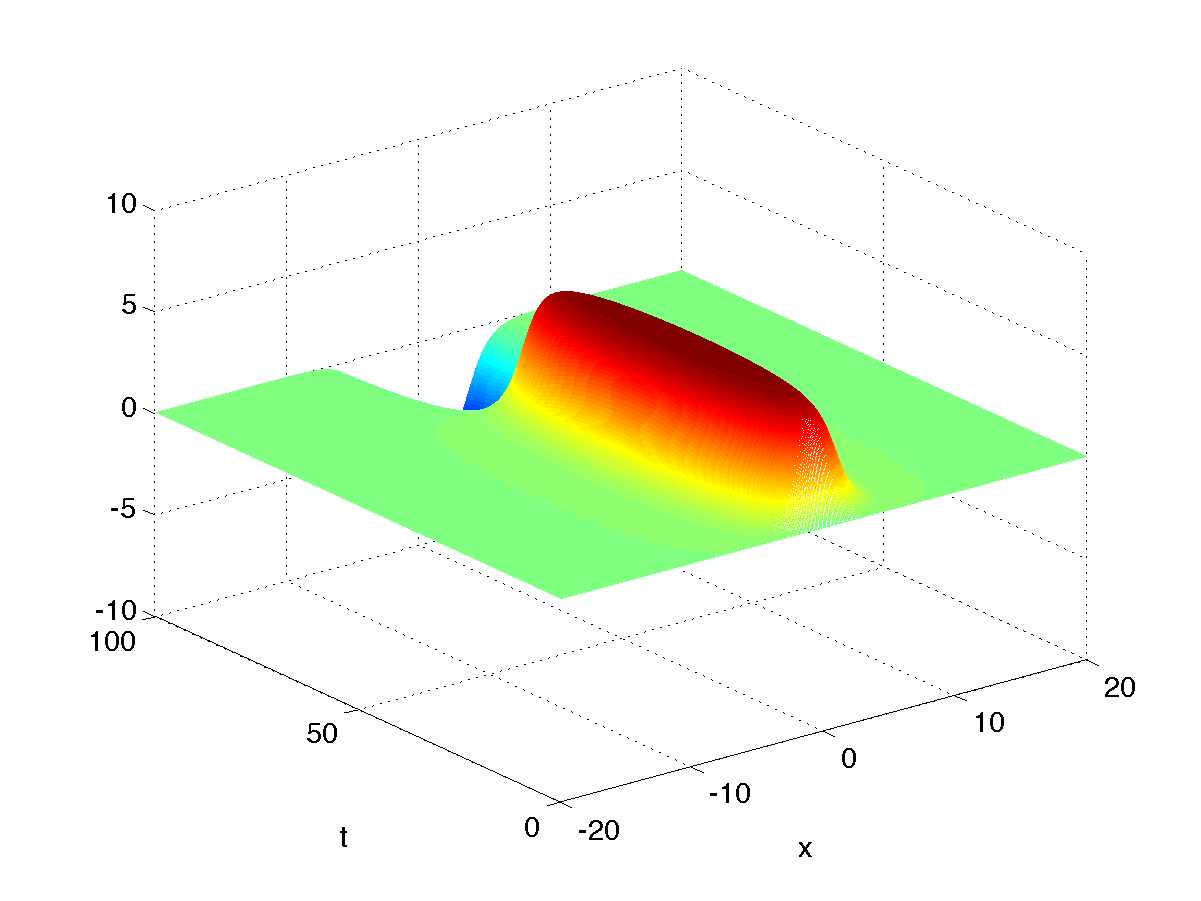}}
\caption{Dirichlet boundary conditions. Difference with the initial value for the numerical Hamiltonian and augmented Hamiltonian (left plot) when solving problem (\ref{sineG})-(\ref{sineG0}) with $\gamma=1$, by using HBVM(1,1) with stepsize $h=0.1$, along with the computed solution (right plot).}
\label{test2}
\end{figure}

\begin{figure}[t]
\centerline{\includegraphics[width=9cm,height=7.5cm]{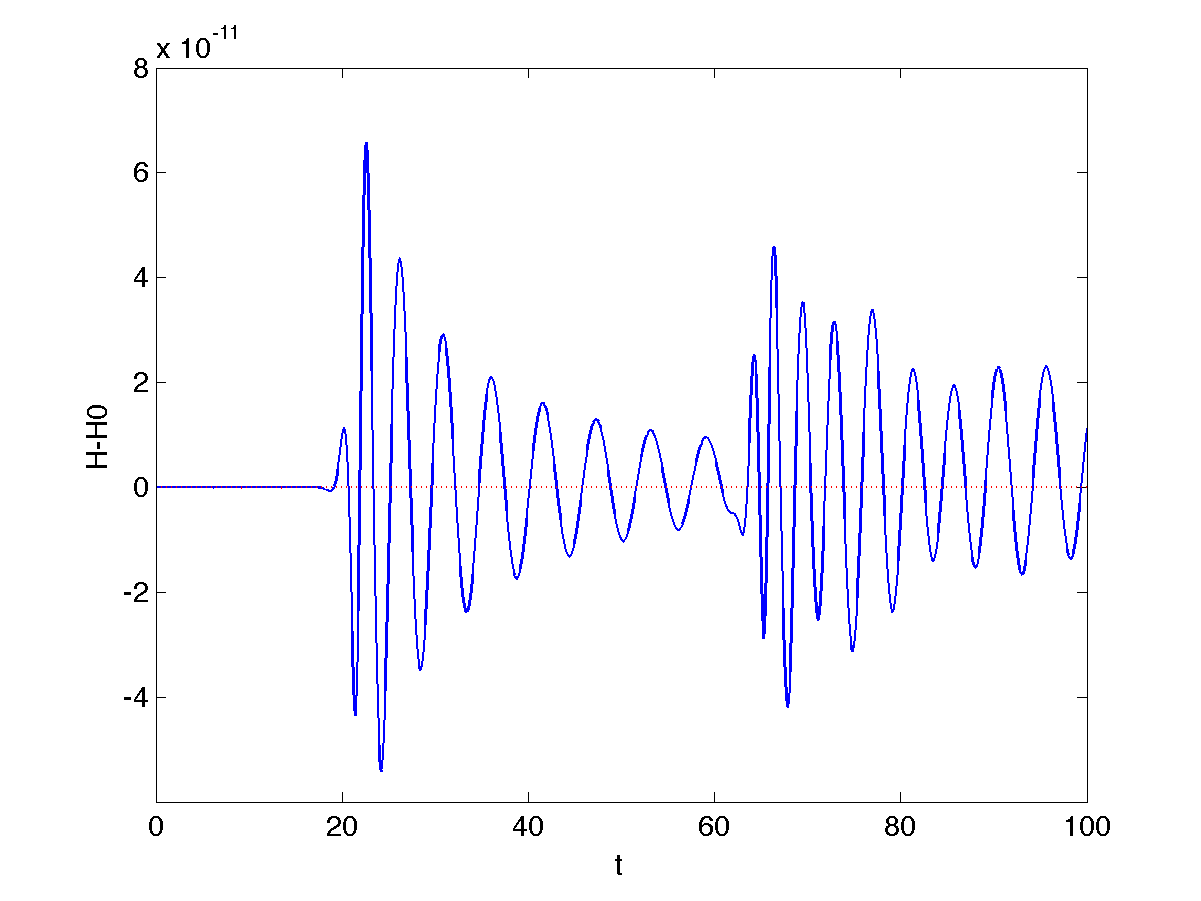}\quad\includegraphics[width=9cm,height=7.5cm]{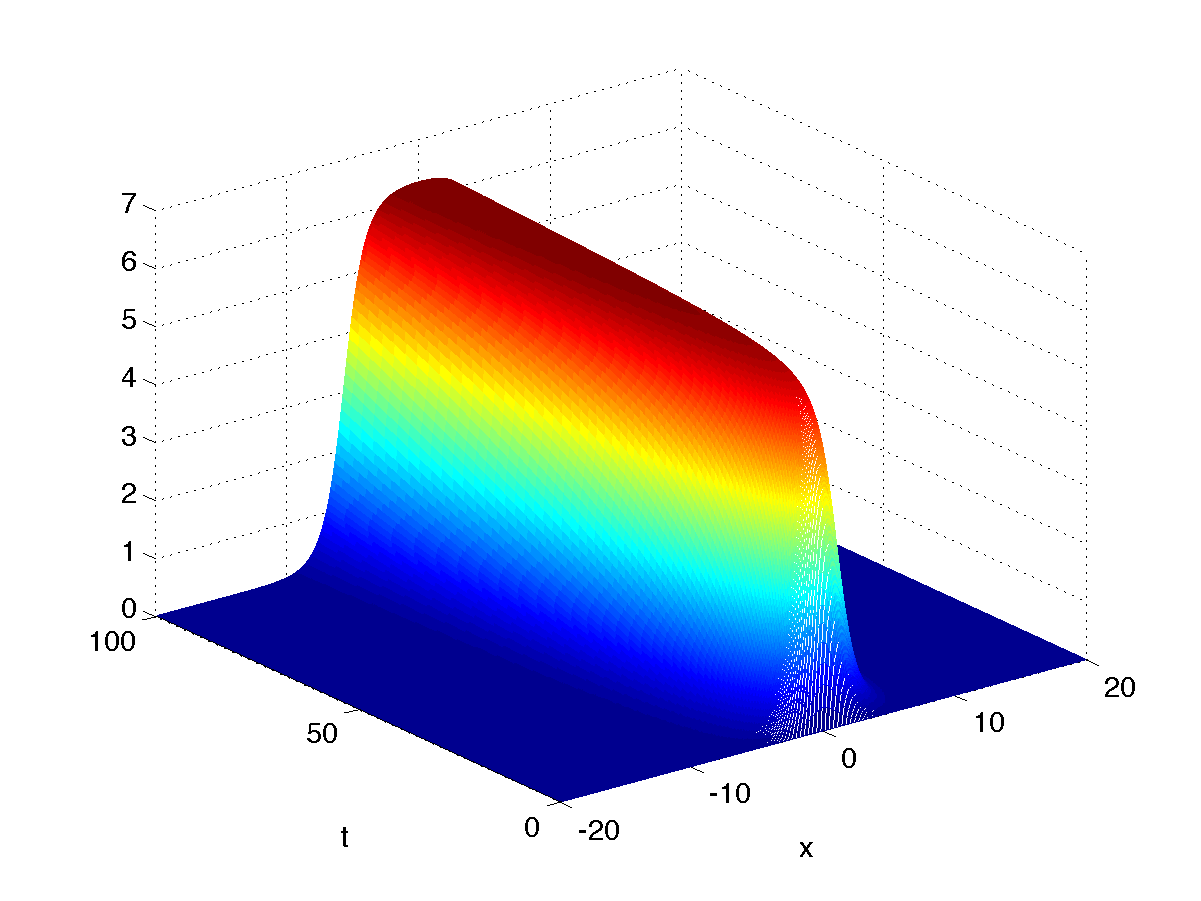}}
\caption{Dirichlet boundary conditions. Difference with the initial value for the numerical Hamiltonian and augmented Hamiltonian (left plot) when solving problem (\ref{sineG})-(\ref{sineG0}) with $\gamma=1$, by using HBVM(5,1) with stepsize $h=0.1$, along with the computed solution (right plot).}
\label{test2_1}
\end{figure}

Completely similar results are obtained by using the same methods (and with the same stepsize $h$), when Dirichlet boundary conditions are prescribed for (\ref{sineG})-(\ref{sineG0}):
\begin{itemize}

\item on the left of Figure~\ref{test2}, there is the plot of $H(\bq_n,\bp_n,t_n)-H(\bq_0,\bp_0,0)$ (see  (\ref{Hassdv})) and $\tH(\bq_n,\bp_n,\tilde{q}_n,\tilde{p}_n)-\tH(\bq_n,\bp_n,\tilde{q}_0,\tilde{p}_0)$ (see (\ref{Hmod})), when using the HBVM(1,1) method. Both differences are quite large and almost overlapping. As a result, the computed numerical solution, shown in the right plot of  Figure~\ref{test2}, is wrong;

\item on the left of Figure~\ref{test2_1}, there is the plot of $H(\bq_n,\bp_n,t_n)-H(\bq_0,\bp_0,0)$ (see  (\ref{Hassdv})) and $\tH(\bq_n,\bp_n,\tilde{q}_n,\tilde{p}_n)-\tH(\bq_n,\bp_n,\tilde{q}_0,\tilde{p}_0)$ (see (\ref{Hmod})), when using the HBVM(5,1) method. The augmented Hamiltonian (\ref{Hmod}) is now conserved, whereas the original Hamiltonian (\ref{Hassdv}) undergoes  small oscillations around its initial value. The computed solution, shown in the right plot of Figure~\ref{test2_1}, is now correct.
\end{itemize}

\begin{figure}[t]
\centerline{\includegraphics[width=9cm,height=7.5cm]{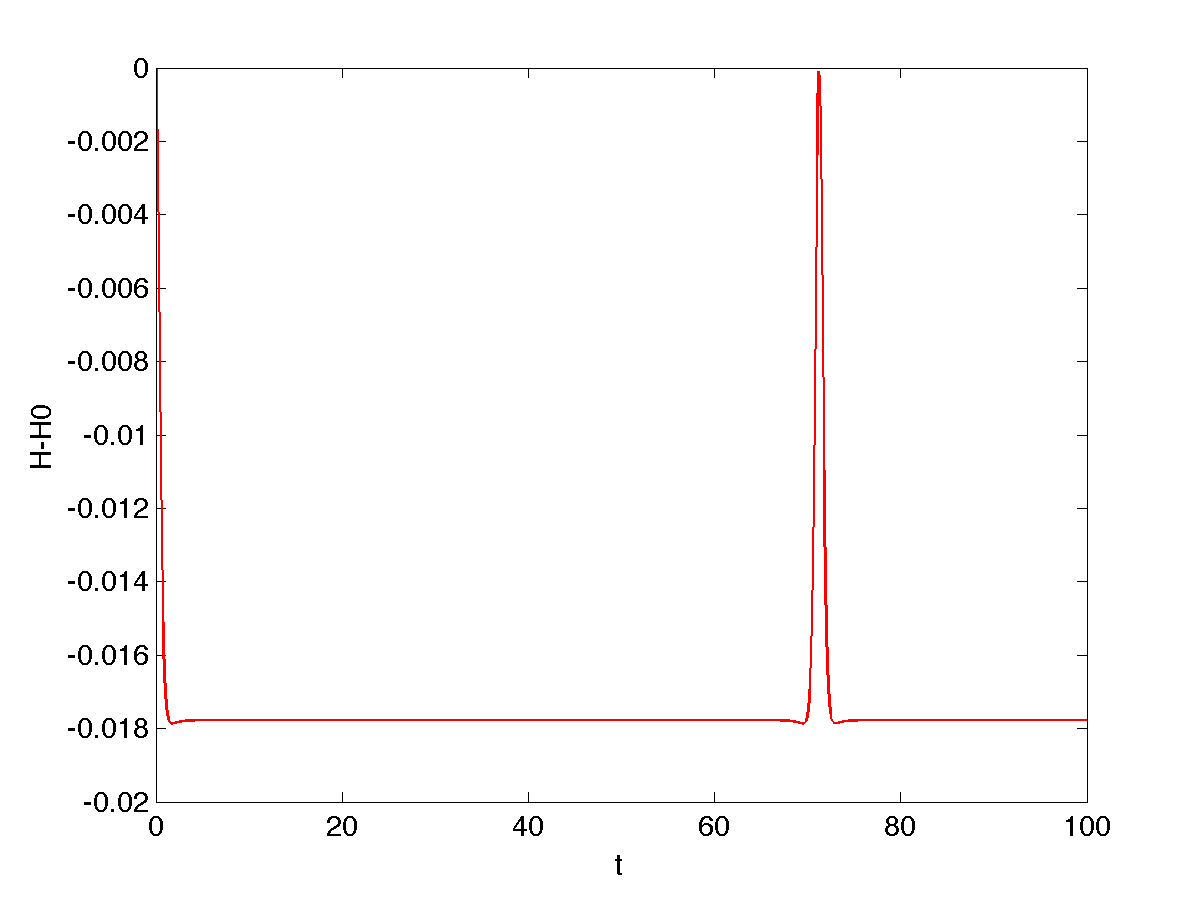}\quad\includegraphics[width=9cm,height=7.5cm]{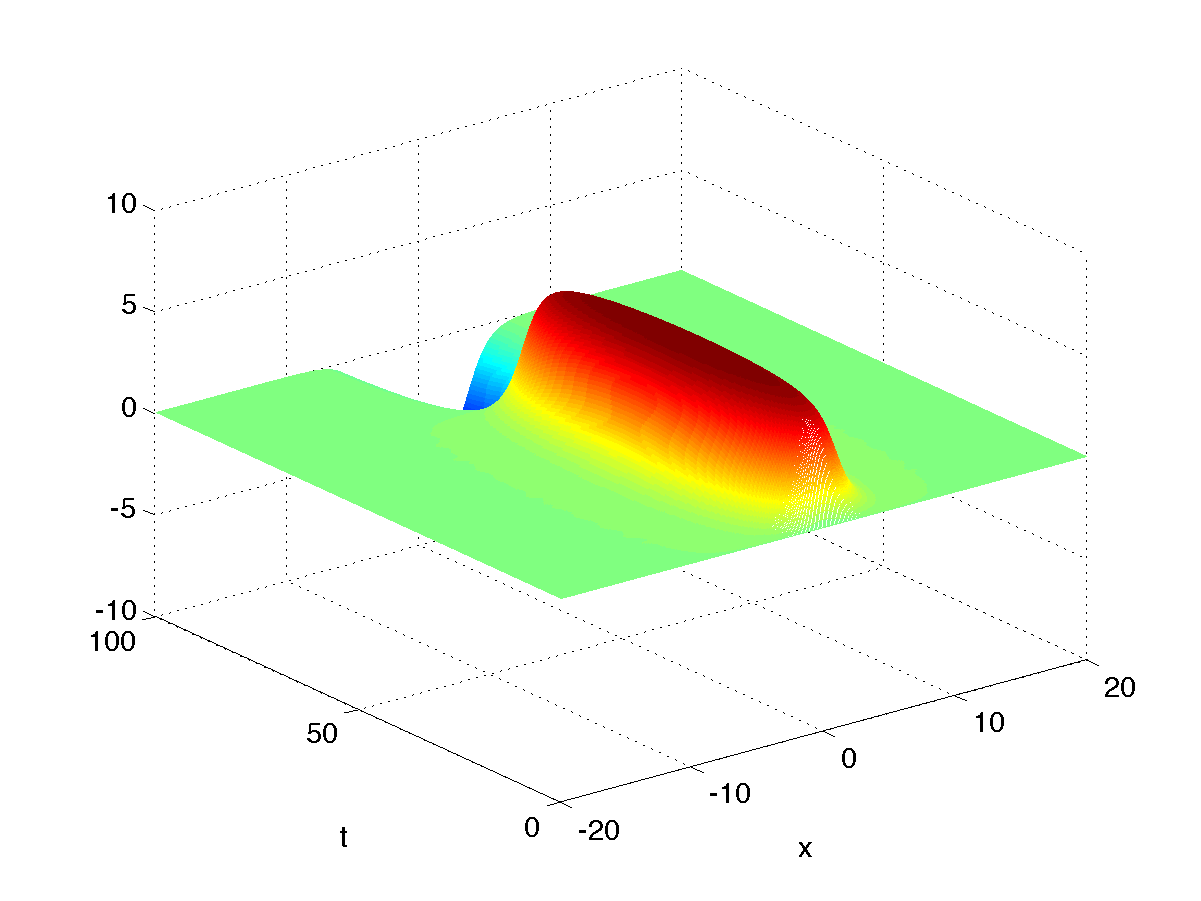}}
\caption{Neumann boundary conditions. Difference with the initial value for the numerical Hamiltonian and augmented Hamiltonian (left plot) when solving problem (\ref{sineG})-(\ref{sineG0}) with $\gamma=1$, by using HBVM(1,1) with stepsize $h=0.1$, along with the computed solution (right plot).}
\label{test3}
\end{figure}

\begin{figure}[t]
\centerline{\includegraphics[width=9cm,height=7.5cm]{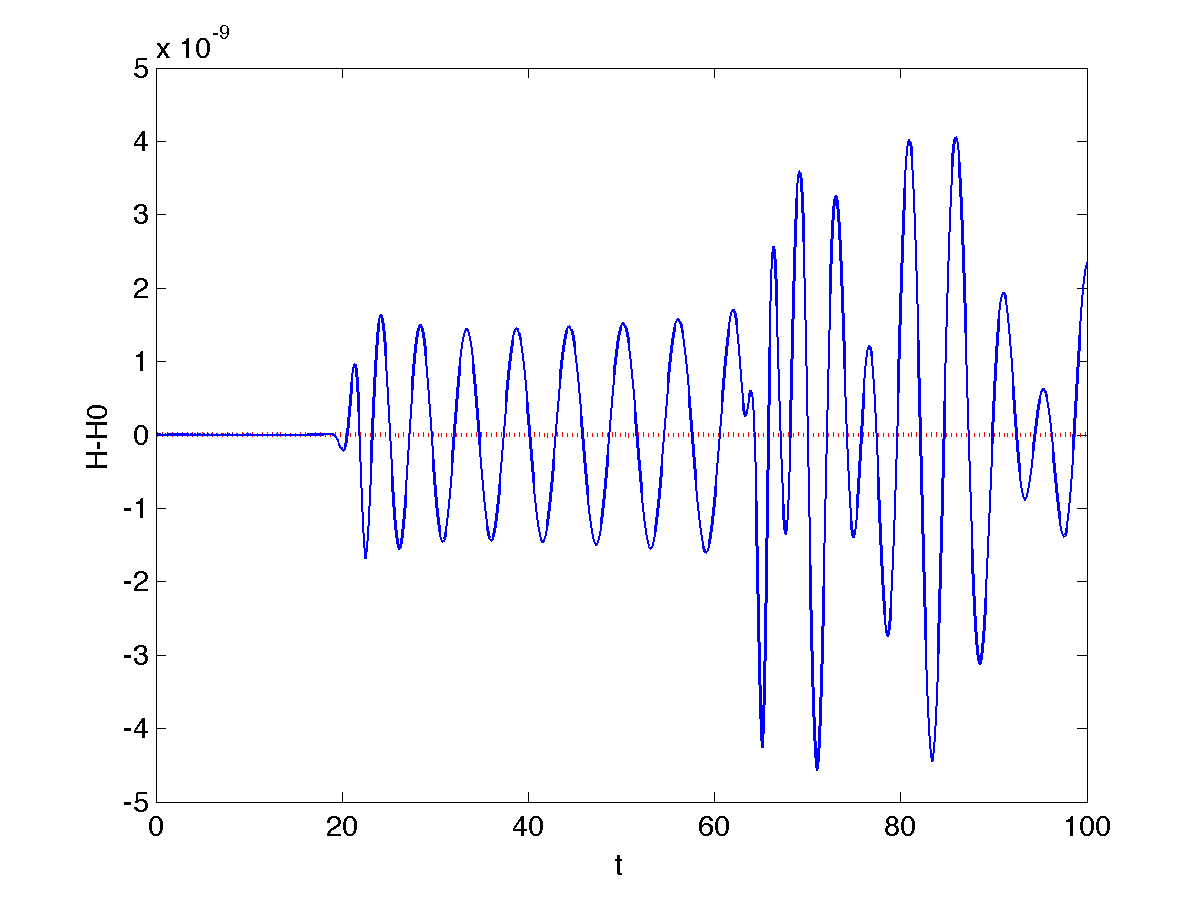}\quad\includegraphics[width=9cm,height=7.5cm]{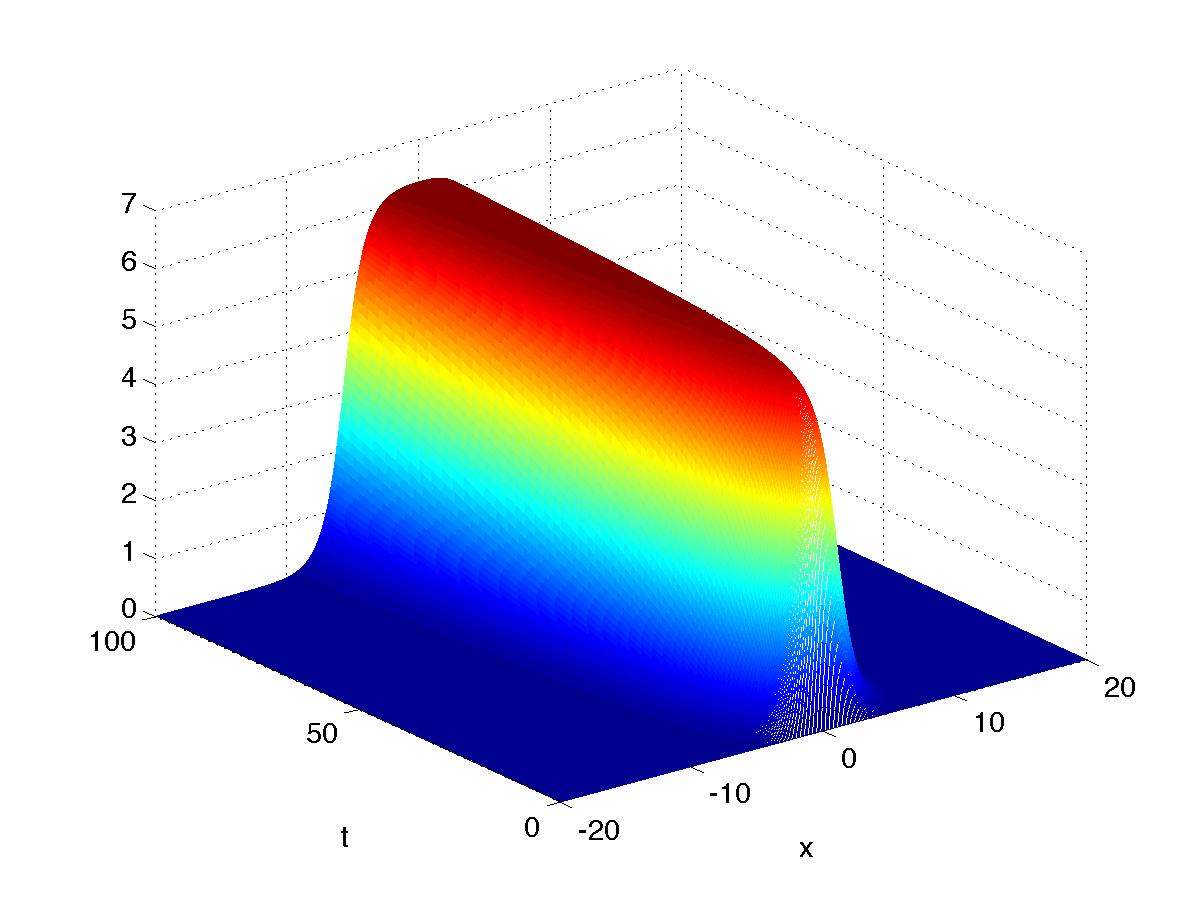}}
\caption{Neumann boundary conditions. Difference with the initial value for the numerical Hamiltonian and augmented Hamiltonian (left plot) when solving problem (\ref{sineG})-(\ref{sineG0}) with $\gamma=1$, by using HBVM(5,1) with stepsize $h=0.1$, along with the computed solution (right plot).}
\label{test3_1}
\end{figure}

Analogous results are obtained when Neumann boundary conditions are prescribed for (\ref{sineG})-(\ref{sineG0}). In fact, by considering the same methods and stepsize $h$:
\begin{itemize}

\item on the left of Figure~\ref{test3}, there is the plot of $H(\bq_n,\bp_n,t_n)-H(\bq_0,\bp_0,0)$ (see  (\ref{Hassnv})) and $\tH(\bq_n,\bp_n,\tilde{q}_n,\tilde{p}_n)-\tH(\bq_n,\bp_n,\tilde{q}_0,\tilde{p}_0)$ (see  (\ref{Hmod1})), when using the HBVM(1,1) method. Both of them are quite large and almost overlapping. As a result, the computed numerical solution, shown in the right plot of Figure~\ref{test3}, is wrong;

\item on the left of Figure~\ref{test3_1}, there is the plot of $H(\bq_n,\bp_n,t_n)-H(\bq_0,\bp_0,0)$ (see  (\ref{Hassnv})) and $\tH(\bq_n,\bp_n,\tilde{q}_n,\tilde{p}_n)-\tH(\bq_n,\bp_n,\tilde{q}_0,\tilde{p}_0)$ (see  (\ref{Hmod1})), when using the HBVM(5,1) method. The augmented Hamiltonian (\ref{Hmod1}) is now conserved, whereas the original Hamiltonian (\ref{Hassnv}) undergoes  small oscillations around its initial value. The computed solution, shown in the right plot of Figure~\ref{test3_1}, is now correct.
\end{itemize}

We now highlight the potentialities of the Fourier-Galerkin space approximation, with respect to the finite-difference one, when periodic boundary conditions are prescribed for the problem: in fact, the Fourier approximation (\ref{H2N}) to the Hamiltonian converges more than exponentially in the number $N$ of Fourier modes, whereas the finite-difference approximation (\ref{Hdisp}) converges only quadratically in $\Dx$. Since also HBVM(5,1) is second order, we then compare the use of such a method, with stepsize $h=40/\ell$ in time and for a total of $\ell$ time-steps, for solving problem (\ref{sineG})-(\ref{sineG0}), with $\gamma=1$ and periodic boundary conditions, by using:
\begin{itemize}
\item the second-order finite-difference spatial discretization with $\ell$ mesh points (with this choice, one has $\Dx=h$);
\item the Fourier-Galerkin approximation with $N=100$, and $m=200$ spatial grid-points, which we maintain fixed independently of the choice of $\ell$. This because the obtained spatial approximation yields a far more accurate approximation than the one corresponding to the time discretization.
\end{itemize}
Table~\ref{errors} summarizes the obtained results:  both methods are globally second-order accurate, even though the values of $N$ and $m$ are kept fixed in the second case (thus confirming the well known exponential convergence of the Fourier approximation). Moreover, by comparing the maximum error in the finite-difference case (FD-error) and in the Fourier-Galerkin approach (FG-error), one sees that the latter is much more favourable than the former.
\begin{table}[t]
\caption{Comparing finite-difference (FD) and Fourier-Galerkin (FG) errors.}
\label{errors}
\small
\centerline{\begin{tabular}{|r|rr|rr|}
\hline
$\ell$       & FD-error & rate & FG-error & rate \\
\hline
     400  & 1.4486e-01  &   --   & 1.7883e-03 &     -- \\
     800  & 3.6900e-02  & 1.97  & 4.4985e-04 &  1.99 \\
   1600  & 9.2702e-03  & 1.99  & 1.1262e-04 &  2.00 \\
   3200  & 2.3204e-03  & 2.00  & 2.8171e-05 &  2.00 \\
\hline   
\end{tabular}}
\end{table}

This fact, allows us to perform a further numerical experiment, where we compare some (practically) energy-conserving HBVMs, 
with well known explicit methods of the same order, for solving problem (\ref{sineG})-(\ref{sineG0}), with $\gamma=1$ and periodic boundary conditions, on the time interval $[0,100]$. In more details, we compare the following methods:
\begin{description}
\item[order 2:] the (practically) energy-conserving HBVM(5,1) method, and the symplectic St\"ormer-Verlet method (SV2);

\item[order 4:] the (practically) energy-conserving HBVM(6,2) method, and the composition method (SV4) based on the symplectic St\"ormer-Verlet method (each step requiring 3 steps of the basic method), according to \cite[page\,44]{HLW06};

\item[order 6:] the (practically) energy-conserving HBVM(9,3) method, and the composition method (SV6) based on the symplectic St\"ormer-Verlet method (each step requiring 9 steps of the basic method), according to \cite[page\,44]{HLW06}.
\end{description}
To compare the methods, we construct a corresponding {\em Work-Precision Diagram}, by following the standard used in the {\em Test Set for IVP Solvers} \cite{testset}. In more details, we plot the accuracy, measured in terms of the maximum absolute error, w.r.t. the execution time. All tests have been done by using Matlab v.\,2014b, running on a dual core i7 at 2.8 GHz computer with 8GB of central memory. The curve of each method is obtained by using $k$ (logarithmically) equispaced steps between $h_{\min}$ and $h_{\max}$, as specified in  Table~\ref{wpdpar}.\footnote{Larger values of $h_{\max}$ for the explicit methods (see Table~\ref{wpdpar}) are not allowed because of stability reasons.} When the stepsize used does not exactly divide the final time $T=100$, the nearest mesh-point is considered.
\begin{table}[t]
\caption{Parameters used for constructing Figures~\ref{wpd1} and \ref{wpd2}.}
\label{wpdpar}
\centerline{\begin{tabular}{|l|rrr|}
\hline
Method       & $h_{\max}$ & $h_{\min}$ & $k$ \\
\hline
HBVM(5,1) & 0.5 & 0.003 & 10\\
HBVM(6,2) & 0.5 & 0.1     & 4\\
HBVM(9,3) & 1    & 0.25   & 4\\
SV2           & 0.1 & 0.0006 & 13\\
SV4           & 0.1 & 0.007   & 7\\
SV6           & 0.1 & 0.01     & 5\\
\hline   
\end{tabular}}
\end{table}

\begin{figure}[t]
\centerline{\includegraphics[width=18cm,height=15cm]{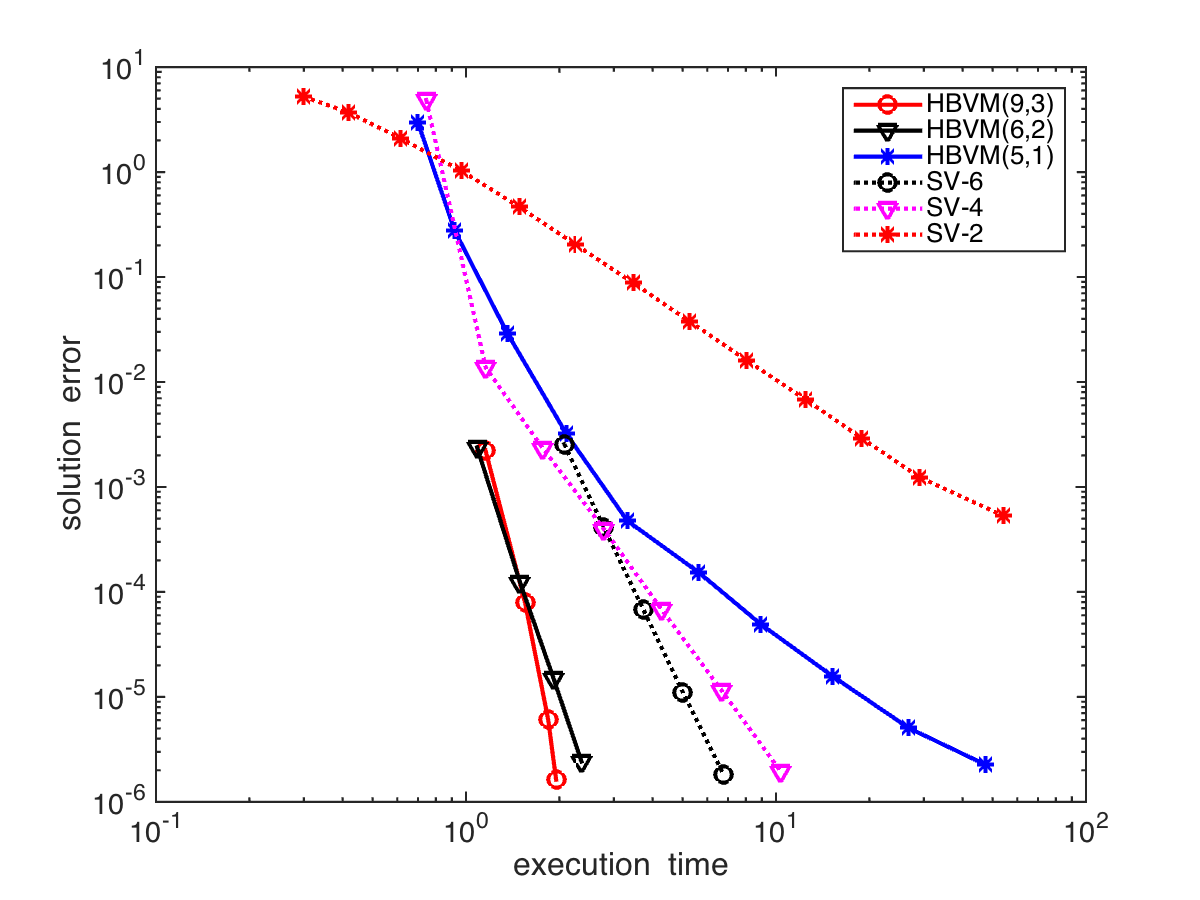}}
\caption{{\em Work-Precision Diagram} ~for problem (\ref{sineG})-(\ref{sineG0}).}
\label{wpd1}
\end{figure}

\begin{figure}[t]
\centerline{\includegraphics[width=18cm,height=15cm]{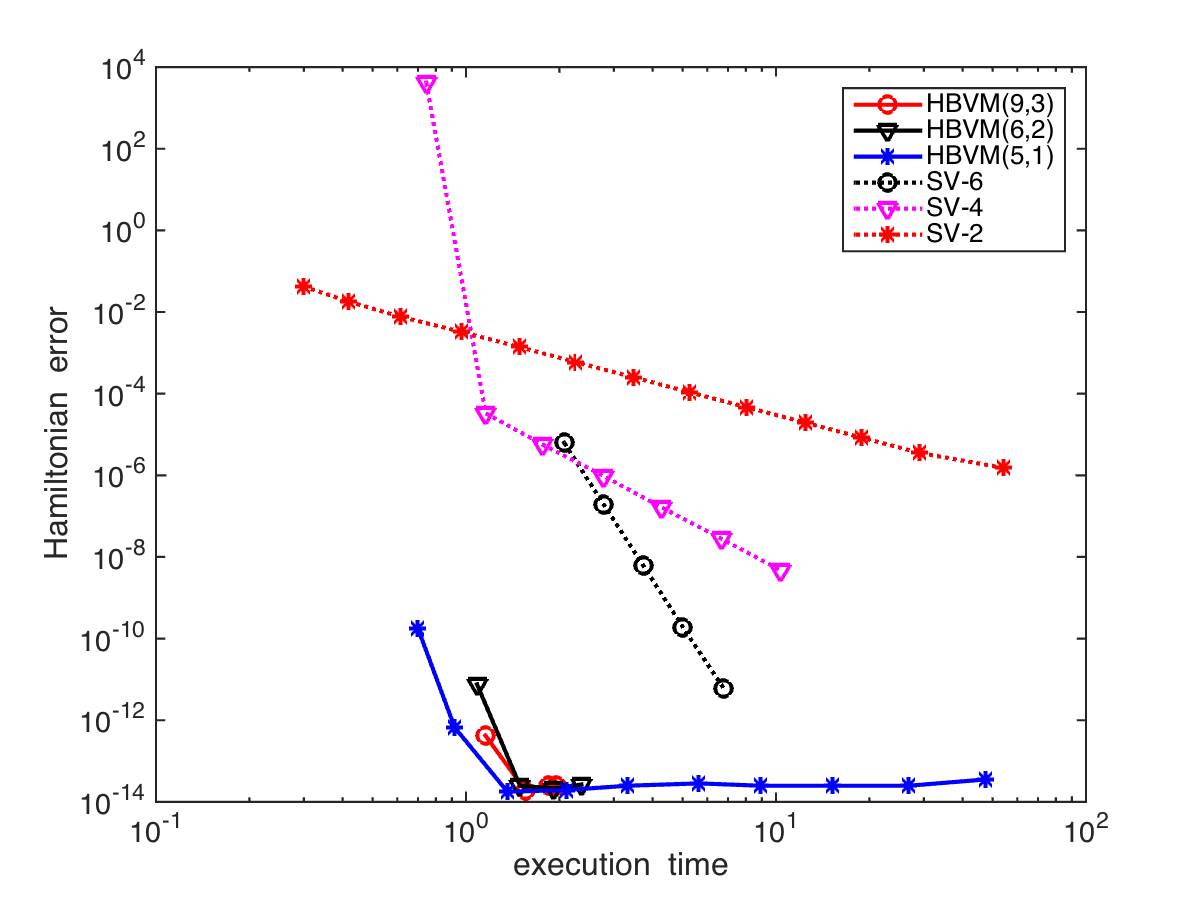}}
\caption{Hamiltonian error versus execution time for problem (\ref{sineG})-(\ref{sineG0}).}
\label{wpd2}
\end{figure}

Figure~\ref{wpd1} summarizes the obtained results, and one sees that the (practically) energy-conserving HBVMs are competitive, even w.r.t. explicit solvers of the same order. For sake of completeness, in Figure~\ref{wpd2}, we plot the corresponding Hamiltonian error versus the execution time, thus confirming that HBVMs are practically energy conserving also for non polynomial Hamiltonians: in fact, taking aside the coarser  time steps, all methods have a Hamiltonian error which is within roundoff errors. On the contrary, for the other methods the decrease of the Hamiltonian error matches their order.

\section{Conclusions}\label{final}

In this paper, we have compared the conservation properties of the semilinear wave equation with the corresponding ones obtained after semi-discretization of the space variable, both when considering a finite-difference and a spectral space discretization. 
When a finite-difference space discretization is considered, we have also studied the case when non-periodic boundary conditions are prescribed for the problem.

The conservation properties of the semi-discrete problem can be conveniently inherited by the numerical solution provided by energy-conserving methods in the HBVMs class. Such methods turn out to be computationally appealing, since they result to be competitive even w.r.t. to explicit methods, and allow a safer approximation of the solution, when energy conservation is an issue, as is confirmed by a few numerical tests on the sine-Gordon equation with a soliton-like solution.

The arguments can be extended in a quite straightforward way to other Hamiltonian partial differential equations, e.g., the Schr\"odinger equation  (as is sketched in the Appendix),  which will be the subject of future investigations. Also a more comprehensive study of Fourier-Galerkin space semi-discretization, when non periodic boundary conditions are prescribed, will be considered in future investigations. 

A further direction of investigation will concern the conservation of multiple invariants for the semi-discrete problem, by means of arguments similar to those used in \cite{BI12,BS13}.

\section*{Appendix}
We here sketch the basic facts that allow an extension of the analysis carried out for the semilinear wave equation (\ref{wave}), to different Hamiltonian PDEs. In particular, we here consider the nonlinear Schr\"odinger equation (in dimensionless form),
\begin{equation}\label{nls}
\ii\psi_t + \psi_{xx} + 2\kappa |\psi|^2\psi = 0, \qquad (x,t)\in (0,1)\times(0,\infty),  \qquad \psi(x,0) \quad \mbox{given},
\end{equation} 
where $\ii$ denotes, as usual, the imaginary unit. By setting $$\psi = u + \ii v,$$ one then obtains the real form of (\ref{nls}),
\begin{eqnarray}\label{nlse}
u_t &=& -v_{xx} -2\kappa (u^2+v^2)v, \qquad (x,t)\in (0,1)\times(0,\infty),\\
v_t &=&  u_{xx}+2\kappa(u^2+v^2)u,  \nonumber
\end{eqnarray}
which is Hamiltonian with Hamiltonian (compare with (\ref{H1}))
\begin{eqnarray}\nonumber
\calH[u,v](t)&=&\frac{1}{2}\int_0^1 \left[u_x^2(x,t) +v_x^2(x,t) -\kappa\left(u^2(x,t) +v^2(x,t)\right)^2 \right]\dd x \\
&\equiv& \int_0^1 E(x,t)\,\dd x. \label{Hnls}
\end{eqnarray}
In fact, (\ref{nlse}) can be formally recast as in (\ref{Hwave})--(\ref{Hwave2}), with the new Hamiltonian function (\ref{Hnls}).
In order to be able to repeat for (\ref{nlse}) the arguments seen for the Hamiltonian semi-discretization of (\ref{wave}), with either periodic, or Dirichlet, or Neumann boundary conditions, it is enough to derive the conservation laws corresponding to (\ref{EtFx}) and (\ref{newH}). Concerning the former conservation law, from (\ref{Hnls}) and (\ref{nlse}) one obtains:
\begin{eqnarray*}
E_t(x,t) &=& u_x(x,t)u_{xt}(x,t) + v_x(x,t)v_{xt}(x,t) - \overbrace{ 2\kappa\left(u^2(x,t) +v^2(x,t)\right)u(x,t)}^{=\,v_t-u_{xx}} u_t(x,t)\\
&&\underbrace{-2\kappa\left(u^2(x,t) +v^2(x,t)\right)v(x,t)}_{=\,u_t+v_{xx}} v_t(x,t)\\
&=& u_x(x,t)u_{xt}(x,t) +u_t(x,t)u_{xx}(x,t) + v_x(x,t)v_{xt}(x,t) +v_t(x,t)v_{xx}(x,t) \\
&=& \left( u_x(x,t)u_t(x,t) \right)_x + \left( v_x(x,t)v_t(x,t) \right)_x ~\equiv~ -F_x(x,t).
\end{eqnarray*}
Consequently, in place of (\ref{EtFx}) one obtains:
$$E_t(x,t)+F_x(x,t), \qquad F(x,t) = -u_x(x,t)u_t(x,t)-v_x(x,t)v_t(x,t).$$
Similarly, taking into account (\ref{Hnls}) and (\ref{nlse}), the analogous of (\ref{newH}) is given by:
\begin{eqnarray*}
\calH[u,v](t)&=& \frac{1}{2}\int_0^1 \left[u_x^2(x,t) +v_x^2(x,t) -\kappa\left(u^2(x,t) +v^2(x,t)\right)^2 \right]\dd x\\
&=&\frac{1}2\int_0^1 \left[-u(x,t)u_{xx}(x,t) +(u(x,t)u_x(x,t))_x \right.\\
&& \left. -v(x,t)v_{xx}(x,t) +(v(x,t)v_x(x,t))_x -\kappa\left(u^2(x,t) +v^2(x,t)\right)^2 \right]\dd x\\ 
&=&-\frac{1}2\int_0^1 \left[u(x,t)u_{xx}(x,t) +v(x,t)v_{xx}(x,t)  +\kappa\left(u^2(x,t) +v^2(x,t)\right)^2 \right]\dd x\\
&& + \frac{1}{2}\left[ u(1,t)u_x(1,t)-u(0,t)(t)u_x(0,t) +v(1,t)v_x(1,t)-v(0,t)(t)v_x(0,t)\right].
\end{eqnarray*}
The arguments for the Hamiltonian semi-discretization of (\ref{nlse}) can then be repeated, {\em mutatis mutandis}, almost verbatim as seen for (\ref{wave}), both when considering a finite-difference and a Fourier-Galerkin space approximation.

\end{document}